\documentclass[11pt]{article}   
\usepackage[latin1]{inputenc}
\usepackage[T1]{fontenc}
\usepackage{amssymb,theorem,amsmath,amsfonts,amscd,geometry,graphics,epsf,epsfig,psfrag,multicol,path,
changebar, array, delarray, enumerate, mathrsfs, setspace}
\usepackage[all]{xy}
\makeindex

\newcommand\qed{\hfill$\square$}

\newtheorem{The}{Theorem}[section]
\newtheorem{Pro}[The]{Proposition}
\newtheorem{Lem}[The]{Lemma}
\newtheorem{Def}[The]{Definition}
\newtheorem{Cor}[The]{Corollary}

\newtheorem{Rem}[The]{Remark}

\numberwithin{equation}{section}

\DeclareMathOperator{\End}{End}
\DeclareMathOperator{\Lie}{Lie}
\DeclareMathOperator{\Spec}{Spec}
\DeclareMathOperator{\Spf}{Spf}
\DeclareMathOperator{\Int}{Int}
\DeclareMathOperator{\Nilp}{Nilp}

\newcommand\Zp{\mathbb{Z}_p}

\begin{document}
\title{ Antispecial cycles on the Drinfeld upper half plane and degenerate Hirzebruch-Zagier cycles}
\author{Ulrich Terstiege}
\date{}
\maketitle
\begin{abstract} We define the notion of antispecial cycles on the Drinfeld upper half plane  in analogy to the notion of special cycles in \cite{KR1}. We determine equations for antispecial cycles and calculate the intersection multiplicity of two antispecial cycles. The result is applied to calculate  the intersection multiplicity of certain degenerate Hirzebruch-Zagier cycles. Finally we compare this intersection multiplicity to certain representation densities.
 \end{abstract}

\section*{Introduction}
Let $k$ be an algebraically closed field of characteristic $p>2$, and let $W=W(k)$ be its ring of Witt vectors. Let $B$ be a division quaternion algebra over $\mathbb{Q}_p$, and let $O_B$ be its ring of integers. According to Drinfeld (\cite{D}), a \emph{special formal $O_B$-module } over a $W$-scheme $S$ is a $p$-divisible formal group $X$ over $S$ 
of dimension 2 and height 4, with an $O_B$-action $\iota: O_B \rightarrow \End_S(X)$ such that the induced 
$\mathbb{Z}_{p^2}\otimes \mathcal{O}_S$-module $\Lie X$ is, locally on $S$, free of rank 1. 
Let $\widehat{\Omega}$ be the formal model of the Drinfeld upper half plane for $\mathbb{Q}_p$, and let 
 $\mathcal{M}=\widehat{\Omega}\times_{\begin{rm}Spf \end{rm}\ \mathbb{Z}_p}\begin{rm}Spf\end{rm}\ W$. Recall that $\mathcal{M}$ represents the following functor on the category Nilp of $W$-schemes $S$ such that $p$ is locally nilpotent in $\mathcal{O}_S$ (comp. \cite{BC} or \cite{KR1}). Let $\mathbb{X}$ be a fixed special formal $O_B$-module over $\Spec k$.  The functor associates to a scheme $S \in \Nilp$ the  set of isomorphism classes of pairs $(X,\varrho)$, where $X$ is a special formal $O_B$-module over $S$ and where 
\[
\varrho: \mathbb{X}\times_{\Spec k}\overline{S} \rightarrow X\times_S \overline{S}
\]
 is an $O_B$-linear quasi-isogeny of height $0$.
Here, 
$\overline{S}=S\times_{\Spec W} \Spec k$. 

We write 
\[\mathbb{Z}_{p^2}= \mathbb{Z}_{p}[\delta]/(\delta^2-\Delta),
\]
 for some  $\Delta\in \mathbb{Z}_{p}^{\times}$ which is not a square and write 
 \[O_B=\mathbb{Z}_{p^2}[\Pi]/(\Pi^2-p),\  \Pi a=a^{\sigma}\Pi\ (a \in \mathbb{Z}_{p^2}).
 \]

 Following \cite{KR1} we call an element $j\in \End(\mathbb{X})\otimes_{\Zp}\mathbb{Q}_p$ \emph{special} if it commutes with the $O_B$-action and its  trace is $0$.

Let $*$ be a $\mathbb{Q}_p$-linear automorphism of order $2$ of $B$. By the theorem of Skolem-Noether there is some $b_* \in B$ such that $*=\Int(b_*)$, i.e. $x^*=b_*xb_*^{-1}$. The element $b_*$ is unique up to multiplication by some element of $\mathbb{Q}_p^{\times}$. We say that an element $j\in \End(\mathbb{X})\otimes_{\Zp}\mathbb{Q}_p$ is \emph{$*$-special} if $j\iota(a)j^{-1}=\iota(a^*) \ \forall a \in B$, and if $s_j:=\iota(b_*^{-1})j$ has trace $0$.  (This means that $s_j$ is special, since it follows from the first condition that $s_j$ commutes with the $O_B$-action.)

Since $*$ has order $2$, it follows that $b_*^2 \in \mathbb{Q}_p$. We will choose $b_*$ in such a way that $b_*^2$ has valuation $0$ or $1$. Precisely one of these possibilities can be realized. In the first case we call $*$ unramified (since $\mathbb{Q}_p(b_*)$ is an unramified quadratic extension of $\mathbb{Q}_p$ in this case), in the second case we call $*$ ramified (since $\mathbb{Q}_p(b_*)$ is a ramified quadratic extension of $\mathbb{Q}_p$ in this case).
Hence $b_*$ is a unit in $O_B$ if $*$ unramified, and $b_*=\Pi\cdot \varepsilon_*$, where $\varepsilon_* $ is a unit in $O_B$, if $*$ is ramified. In the latter case we write $b_*^2=\eta_* p$ for some $\eta_*\in \Zp^{\times}$. Note that since $\varepsilon_*$ is unique up to multiplication by some element of $\Zp^{\times}$, it follows that the quadratic residue character $\chi(\eta_*)$ is well defined.

Let $V[*]$ be the space of $*$-special endomorphisms. It is a quadratic $\mathbb{Q}_p$-vector space via the quadratic form
\[
q(j)= (ps_j)^2.
\]
(The scaling by $p$ is done to facilitate the comparison with special cycles, cf. Theorem \ref{A}. Note that the ambiguity in the choice of $b_*$ leads to an ambiguity in $s_j$ and hence in $q$. But $q$ is unique up to multiplication by some element of $\Zp^{\times, 2}$.)
The quadratic form $q$ induces a bilinear form $\beta$ on $V[*]$ given by $\beta(j_1,j_2)=q(j_1+j_2)-q(j_1)-q(j_2)$.
If $j\in \End(\mathbb{X})\otimes_{\Zp}\mathbb{Q}_p$,  we define the  cycle $Z(j)$ to be the closed formal subscheme of  $\mathcal{M}$ consisting of all pairs $(X,\varrho)$ such that $\varrho\circ  j \circ \varrho^{-1}$ lifts to an endomorphism of $X$. If $j$ is special, resp. $*$-special, we call $Z(j)$ a special resp. $*$-special cycle.

The content of the first six sections of  \cite{KR1} is the description of special cycles by equations and the determination of the intersection product of two special cycles. Our first aim in this paper is to do the same for $*$-special cycles. For unramified $*$ and a $*$-special $j$ we have $Z(j)=Z(s_j)$, which is a special cycle. For arbitrary unramified $*$ the notion of $*$-special cycles is the same as the notion of special cycles. Hence from now on we assume that $*$ is ramified. In this case we call a $*$-special cycle also an \emph{antispecial cycle}. 
The notion of antispecial cycles does not depend on the particular $*$, i.e. all ramified $*$ induce the same  $*$-special cycles.
We now state our main results on antispecial cycles.
\begin{The}\label{A}
Let $p>3$, and let $j\in V[*]$ with $q(j)\neq 0$ and $Z(j)\neq \emptyset$. Then $Z(j)$ is a divisor in $\mathcal{M}$. We have
\[ 
Z(j)=Z(ps_j)^{\rm pure},
\]
 where the upper index ${\rm pure}$ means the associated divisor of $Z(ps_j)$, i.e. the subscheme defined by the ideal sheaf of local sections with finite support. If $j^2 \equiv 0 \mod p$ then this  statement is also true in case $p=3$.
\end{The}

For the proof we show first by considering Dieudonn\a'e modules that $Z(j)$ and $Z(ps_j)$ have the same $k$-valued points.  Then we use the known equations for $Z(ps_j)$ on the one hand and Grothendieck-Messing theory on the other hand to give equations for $Z(j)$.

We define the intersection product $(Z(j_1),Z(j_2))$ of two antispecial cycles $Z(j_1)$ and $Z(j_2)$ using the general definition in \cite{KR1}.  
  
\begin{The}\label{B}
Let $j_1,j_2\in V[*]$  such that $j_1,j_2$ span a $2$-dimensional nondegenerate quadratic $\mathbb{Z}_p$-submodule of $V[*]$. Let 
\[
T:=
\begin{pmatrix}
 q(j_1)& \frac{1}{2}\beta(j_1,j_2) \\ 
 \frac{1}{2}\beta(j_2,j_1) &q(j_2)
\end{pmatrix}.
\]
We suppose that $T$ is $GL_2(\mathbb{Z}_p)$-equivalent to $\text{diag}(\eta_1p^{\beta_1},\eta_2p^{\beta_2})$,
where $\eta_1, \eta_2 \in 
\mathbb{Z}_p^\times$ and $1\leq \beta_1 \leq \beta_2$ for $p>3$, resp. $1<\beta_1 \leq \beta_2$ for $p=3$.
Define $\varepsilon_i\in \mathbb{Z}_p^\times $ and $\alpha_i \in \mathbb{N}$ by 
 $\eta_*\eta_i p^{\beta_i -1}=
\varepsilon_i p^{\alpha_i}$. Then  
\[
\begin{split}
&(Z(j_1),Z(j_2))= \\& \alpha_1+\alpha_2+3 
 -
\begin{cases}
p^{(\alpha_1 +1)/2}+2\frac{p^{(\alpha_1 +1)/2}-1}{p-1} \ & \text{if $\alpha_1$ is odd and 
$\chi(\eta_*\varepsilon_1)=-1$}\\
(\alpha_2-\alpha_1+1)p^{(\alpha_1 +1)/2}+ 2\frac{p^{(\alpha_1 +1)/2}-1}{p-1}& \text{if $\alpha_1$ is odd and 
$\chi(\eta_*\varepsilon_1)=1$}\\
2\frac{p^{\alpha_1/2+1}-1}{p-1} \ & \text{if $\alpha_1$ is even,}
\end{cases}
\end{split}
\]
where $\chi$ denotes the quadratic residue character on $\mathbb{Z}_p^\times$ resp.
$\mathbb{F}_p^\times$.
\end{The}
(Note that by the earlier remarks  the values for $\beta_i$, $\alpha_i$ and $\chi(\eta_*\varepsilon_i)$ do not depend on the choice of $b_*$.)
This theorem is in fact a simple consequence of Theorem \ref{A} and  the formula for the intersection product of special cycles given in \cite{KR1} which only depends on the divisors associated to the special cycles.

Our second aim in this paper is to apply  Theorem \ref{B} to compute the intersection product of certain degenerate intersections of arithmetic Hirzebruch-Zagier cycles. We consider the following moduli problem.
Fix a supersingular formal 
$p$-divisible group $\mathcal{A}$ over $k$ of height 4 and dimension 2 which is equipped with an action 
\[
\iota_0: \mathbb{Z}_{p^2}\rightarrow \begin{rm}End\end{rm}(\mathcal{A}),
\]
such that $\mathcal{A}$ is special with respect to this $\mathbb{Z}_{p^2}$-action.
 We further suppose that $\mathcal{A}$ is equipped with a  polarization
\begin{equation*}
\begin{CD}
\lambda: \mathcal{A} @ >\sim>> \hat{\mathcal{A}},
\end{CD}
\end{equation*}

such that for the Rosati involution $\iota_0(a)^{*}=\iota_0(a)$.
We consider the following functor $\mathcal{M}^{HB}$ on the category Nilp.
It associates to a scheme $S \in \Nilp$  the set of isomorphism classes of the following data.

\begin{enumerate}[(1)]
\item A $p$-divisible group $X$ over $S$, with an action
\[
\iota_0: \mathbb{Z}_{p^2}\rightarrow \begin{rm}End\end{rm}(X),
\]
such that $X$ is special with respect to this $\mathbb{Z}_{p^2}$-action.
\item 
A quasi-isogeny of height zero
\[
\varrho: \mathcal{A}\times_{\Spec k}\overline{S} \rightarrow X\times_S \overline{S},
\]
which commutes with the action of $\mathbb{Z}_{p^2}$   such that the following condition holds.
Let   ${\lambda}_{\overline{S}}: \mathcal{A}_{\overline{S}} \rightarrow
 \hat{\mathcal{A}}_{\overline{S}}$ be the map induced by $\lambda$.
Then we require the existence of an isomorphism $\tilde{\lambda}
:X \rightarrow \hat{X}$
such that for the induced map 
  $\tilde{\lambda}_{\overline{S}}:X_{\overline{S}} \rightarrow \hat{X}_{\overline{S}}$ we have 
${\lambda}_{\overline{S}}=\hat{\varrho}\circ \tilde{\lambda}_{\overline{S}}\circ \varrho$.

\end{enumerate}
This functor is representable by a formal scheme which we also 
call ${\mathcal{M}}^{HB}$. (We note that by \cite{RZ} the scheme ${\mathcal{M}}^{HB}$ can be used to uniformize the completion along the supersingular locus of the Hilbert-Blumenthal moduli surface at an inert prime of the real quadratic field.)
On the isocrystal $N$ of $\mathcal{A}$ we have a perfect symplectic form. We define in this context the space of special endomorphisms 
\[
V^{'}=\{j \in \End(N); j\iota_0(a)=\iota_0(a^{\sigma})j \ \text{ and } \ j^*=j\},
\]
 where $*$ denotes the adjoint with respect to the symplectic form. Then 
$V^{'}$ is a $4$-dimensional vector space over
$\mathbb{Q}_p$ with quadratic form 
\[
Q(j)=j^2.
\]
(Compare \cite{KR2}, §5.)
For a special endomorphism $j\in V^{'}$  we define the special cycle $Z(j)$ as a closed formal subscheme of $\mathcal{M}^{HB}$ as above. In analogy to Theorem \ref{A} we have the following 
\begin{Pro}
  Let $j \in V^{'}$
 be such that $Q(j)\neq0$ and $Z(j)\neq \emptyset$. 
 Then
   $Z(j)$ is a divisor in ${\mathcal{M}}^{HB}$.

\end{Pro}

 Now we fix a special endomorphism $j_1\in V^{'}$ with $j_1^2=\varepsilon_1p$ for some $\varepsilon_1\in \mathbb{Z}_p^{\times}$. Let
\[
V^{'}[j_1]=\{j\in V^{'}\  | \ j \perp j_1 \text{ with respect to the bilinear form associated to }  Q \}.
\] 
Let $j_2,j_3\in V^{'}[j_1]$ be such that the $\mathbb{Z}_p$-span $\begin{bf} j \end{bf}=
\mathbb{Z}_pj_2 + \mathbb{Z}_pj_3$ has rank $2$ as a submodule of $V^{'}$, and such that $Q$ induces a nondegenerate bilinear form on $\bf{j}$. 
We further suppose
that the matrix of the bilinear form on $\bf{j}$ associated
to $Q$ with respect to the basis $  j_2, \ j_3$ is equivalent to 
 $\begin{rm}  diag \end{rm}( \varepsilon_2 p^{\beta_2}, \varepsilon_3 p^{\beta_2})$ with $\varepsilon_i\in \Zp^{\times}$ and $1\leq \beta_2\leq \beta_3$.
We define the intersection product by the \emph{Euler-Poincar\a'e} characteristic of the derived tensor product,
\[
(Z(j_1),Z(j_2),Z(j_3)):= \chi(\mathcal{M}^{HB},\mathcal{O}_{Z(j_1)}\otimes^{\mathbb{L}}\mathcal{O}_{Z(j_2)}
\otimes^{\mathbb{L}}\mathcal{O}_{Z(j_3)}),
\]
which is well defined since $Z(j_1)\cap Z(j_2)\cap Z(j_3)$ is proper over $\Spec k$. 
\begin{Pro}\label{HB}
$Z(j_1)$ is isomorphic to $\mathcal{M}$, such that $Z(j_i)\cap \mathcal{M}$ for $i=2,3$ can be identified with the cycle associated to a  $*$-special endomorphism for $*=\Int (\delta j_1)$ (identifying $B$ with  $\mathbb{Q}_{p^2}[j_1]$). Furthermore 
\[
(Z(j_1),Z(j_2),Z(j_3))=(Z(j_2)\cap \mathcal{M},Z(j_3)\cap \mathcal{M}).
\]
 \emph{The latter is given explicitly by  Theorem \ref{B} in case $*=\Int (\delta j_1)$. }
\end{Pro}

Our third aim is to compare the intersection multiplicity $(Z(j_1),Z(j_2),Z(j_3))$ to certain representation densities.
To formulate the result,
recall that, for $S\in \begin{rm}  Sym \end{rm}_m(\mathbb{Z}_p)$ and 
$T\in \begin{rm}  Sym \end{rm}_n(\mathbb{Z}_p)$
with $\begin{rm}  det \end{rm}(S)\neq 0$ and $\begin{rm}  det \end{rm}(T)\neq 0$,
 the representation density is defined as 
\[
\alpha_p(S,T)= 
\operatorname*{lim}_{t\rightarrow\infty} p^{-tn(2m-n-1)/2} \mid \{x \in M_{m,n}(\mathbb{Z}/p^t\mathbb{Z});
\ S[x]-T \in p^t\begin{rm}  Sym \end{rm}_n(\mathbb{Z}_p)\}\mid.
\]
For $S$ as above, let 
\[
S_r=
\begin{pmatrix}
 S \\ 
 & 1_r \\
 & & -1_r
\end{pmatrix}.
\]
Then there is a rational function $A_{S,T}(X)$ of $X$ such that 
\[
\alpha_p(S_r,T)= A_{S,T}(p^{-r}).
\]
Let 
\[
\alpha_p^{'}(S,T)=\frac{\partial}{\partial X}(A_{S,T}(X)) \arrowvert_{X=1}.
\]

(Compare \cite{KR1}, §7.)
Recall our assumption that the matrix of the bilinear form on $\bf{j}$ associated
to $Q$ with respect to the basis $  j_2, \ j_3$ is equivalent to 
 $\begin{rm}  diag \end{rm}( \varepsilon_2 p^{\beta_2}, \varepsilon_3 p^{\beta_2})$ where $\varepsilon_i \in \mathbb{Z}_p^{\times}$ and $1\leq \beta_2\leq \beta_3$. 
 Now let $T=\begin{rm}  diag \end{rm}( \varepsilon_1 p, \varepsilon_2 p^{\beta_2},  \varepsilon_3 p^{\beta_3})$. (Then $T$ equals the matrix of the bilinear form associated to the quadratic form $Q$ restricted to the $\Zp$-submodule
 ${\bf j^{'}}:=\Zp j_1 \oplus \Zp j_2 \oplus \Zp j_3$ of $V^{'}$ for a suitable basis of $\bf{j^{'}}$ consisting of a suitable basis of $\bf{j}$ together with $j_1$.)
  Let $S=\begin{rm}  diag \end{rm}(1,-1,1,-\Delta)$.
We show the following

\begin{The}\label{C}
\[
(Z(j_1),Z(j_2),Z(j_3))= -\frac{p^4}{(p^2+1)(p^2-1)}\alpha_p^{'}(S,T).
\]
\end{The} 

The proof is by explicit calculation of the r.h.s. and comparing it to the expression given by  Theorem \ref{B} for the l.h.s.. The calculation of the r.h.s.  uses a combination of a lemma of Shimura (\cite{S}) and a formula for $\alpha_p(\tilde{S}_r,T)$ in case $\tilde{S}=\text{diag}(1,-1,1,-1)$ given by Katsurada (\cite{Ka}) which together allow us to calculate $\alpha_p^{'}(S,T)$ for our $S$ and $T$ explicitly. 
\newline

In [KR2] an analogous formula for $(Z(j_1),Z(j_2),Z(j_3))$ is proved in case $Q(j_1)\in \Zp^{\times}$. A conjecture of Kudla and Rapoport states that this formula holds in general provided that
 ${\bf j^{'}}=\Zp j_1  \oplus \Zp j_2 \oplus \Zp j_3$ is of rank $3$ and $Q$ induces a nondegenerate bilinear form on ${\bf j^{'}}$.
Theorem \ref{C} confirms a special case of this conjecture.
\newline 
 
The paper is divided into five sections. The first section introduces some linear algebra concerning Dieudonné modules and the notion of $*$-special endomorphisms. The second section introduces $\mathcal{M}$ and investigates antispecial cycles and their local equations. In the third section we  investigate intersection products of antispecial cycles and prove Theorem \ref{B}. In the fourth section we discuss the application to arithmetic Hirzebruch-Zagier cycles, and in the fifth section we prove Theorem \ref{C}. 
\newline

I want to conclude this introduction by thanking  those people who helped me to write this paper. In particular I thank  M. Rapoport for suggesting this topic and for his stimulating support during the work. My deep thanks go to  U. Görtz for many hours of patient help. Thanks are also due to Prof. Katsurada who gave the crucial hint how to calculate the representation densities and to Prof. Messing for helpful comments on $p$-divisible groups. Finally, I thank Prof. Kudla for alerting me to a mistake in a first version of this paper.
\
%-----------------------------------------------------------------------------------------------------------------------
%-----------------------------------------------------------------------------------------------------------------------

\section{Special Dieudonn\a'e  modules with $O_B$-action, $O_B$-lattices and $*$-special endomorphisms}
As in the introduction, let $k$ be an algebraically closed field of characteristic $p>2$, let $W=W(k)$ be its ring of Witt vectors with
fraction field $W_{\mathbb{Q}}$, and
let $\sigma$ be the Frobenius automorphism of $W$. Also, let $B$ be a quaternion division algebra over
$\mathbb{Q}_p$, and let $O_B$ be its ring of integers, which we identify with
\begin{align*}
O_B = \mathbb{Z}_{p^2}[\Pi]/(\Pi^2-p),\text{} \Pi a = a^{\sigma}\Pi \text{ } \forall a \in \mathbb{Z}_{p^2}.
\end{align*} As in the introduction, we also write
\[
\mathbb{Z}_{p^2}=\mathbb{Z}_p[\delta]/(\delta^2-\Delta) \text{ for some } \Delta \in \mathbb{Z}_p^{\times}
\setminus \mathbb{Z}_p^{\times,2}. 
\]
 We will regard $\mathbb{Z}_{p^2}$ as a subset of $W$ (the set of elements fixed by $\sigma^2$).
Let $(M,F,V)$ be a Dieudonn\a'e module of height 4 and
dimension 2. It is a free $W$-module of rank 4 with a $\sigma$-linear resp. $\sigma^{-1}$-linear endomorphism
$F$ resp. $V$ satisfying $VF=FV=p$, and for which the $k$-vector space $M/VM$ has dimension 2. Now we assume 
that $M$ is
 equipped with an $O_B$-operation, i.e. an action 
$\iota: O_B \rightarrow \begin{rm}End\end{rm}(M)$ commuting with $F$ and $V$. From the action of
$\mathbb{Z}_{p^2}$ we obtain a $\mathbb{Z}/2$-grading,
\[
M=M_0 \oplus M_1,
\] 
where
\[
\begin{split} 
&M_0=\{m\in M \text{ }| \text{ }\iota(a)m=am \text{ }\forall a \in \mathbb{Z}_{p^2}\},  \\
&M_1=\{m\in M \text{ }|\text{ } \iota(a)m=a^\sigma m  \text{ } \forall a \in \mathbb{Z}_{p^2}\}.
\end{split}
\]
If we denote by $N$ the isocrystal of $M$ (i.e. the $W_{\mathbb{Q}}$-vector space
 $M\otimes_{W}W_{\mathbb{Q}}$ equipped with the induced notions für $F$ and $V$ and the $O_B$-action), 
 we also obtain a $\mathbb{Z}/2$-grading,
\[
N=N_0 \oplus N_1.
\] 
%Finally, we require $M$ to be \emph{special}, i.e. the $k$-vector spaces $M_0/VM_1$ and $M_1/VM_0$ are
% of dimension 1. Summarizing we define
 
 \begin{Def}\label{blabla}
 \emph{A} special Dieudonn\a'e  module with $O_B$-action  \emph{is  a 
 Dieudonn\a'e module $(M,F,V)$ over $W$ of height 4 and dimension 2 which is equipped with an $O_B$-action 
 such that the $k$-vector spaces $M_0/VM_1$ and $M_1/VM_0$ are
 one dimensional.}
\end{Def}

 If $(M,F,V)$ is a special 
 Dieudonn\a'e module with $O_B$-action and if $i \in \mathbb{Z}/2$, we say that
 the index $i$ is  $O_B$-\emph{critical} (for $M$), if $VM_i=\iota(\Pi)M_i$. 
  Further $M$ is called $O_B$-\emph{superspecial}, if both indices $0, 1$ are $O_B$-critical. Finally, 
  $M$ is called $O_B$-\emph{ordinary}, if only one index is $O_B$-critical.

\begin{Def} 
\emph{Let $(M,F,V)$ be a special Dieudonn\a'e module with $O_B$-action, and let $N$ be its isocrystal.
 An} $O_B$-lattice in $N$ \emph{is a free $W$-submodule $L$ of rank 4 in $N$, which is spanned by a
 basis of $N$ 
and which is stable under $F$ and $V$ and under $\iota$, and for which 
 the $k$-vector spaces $L_0/{VL_1}$ and $L_1/{VL_0}$ have 
dimension 1. (Here, as usual, $L=L_0\oplus L_1$ is the $\mathbb{Z}/2$- grading obtained from the action of
 $\mathbb{Z}_{p^2}$.)}
\end{Def}

Let us fix $M$ and hence also fix $N$ for this section. Then an $O_B$-lattice $L$ is a special Dieudonn\a'e 
module 
 (with $O_B$-operation) together with an isomorphism of isocrystals
\[
L\otimes_W W_{\mathbb{Q}} \rightarrow N,
\]
 i.e. an isomorphism of vector spaces 
which commutes with the  endomorphisms $F,V$ and the operation of $O_B$.

\begin{Lem}\label{lemmain}
Let $L$ be an $O_B$-lattice, and let $i \in \mathbb{Z}/2$. Then
\begin{enumerate}[(i)]
\item the inclusions 
$pL_i \subset VL_{i+1} \subset L_i$  are both of index 1.
\item the inclusions $pL_{i}\subset\iota(\Pi)L_{i+1} \subset L_{i}$ are both of index 1.
\end{enumerate}
\end{Lem}
\emph{Proof.}  The first statement is essentially the assumption that $L$ is special, for the second
 see  \cite{BC}, chapitre II, Proposition 5.1.
\qed

\begin{Lem}\label{critical}
\begin{enumerate}[(i)]

\item  Any $O_B$-lattice $L$ has an $O_B$-critical index.
\item There exists an $O_B$-superspecial $O_B$-lattice.
\end{enumerate}
\end{Lem}
\emph{Proof.} \emph{i}) See \cite{KR1}, p. 165. 

\emph{ii}) Let $L$ be an $O_B$-lattice, and let $i\in \mathbb{Z}/2$ be $O_B$-critical.
 By the theorem of Dieudonn\a'e (\cite{Z}, Satz 6.26) we can choose  a basis $e_1,e_2$  of $L_i$ satisfying 
$\Pi^{-1}Ve_i=e_i$.  Then the $O_B$-lattice spanned by $e_1,e_2,e_3=Ve_1 \text{ and } e_4= p^{-1}Ve_2$ is
superspecial. \qed

\begin{Def}
\emph{A basis $e_1$, $e_2$, $e_3$, $e_4$ of $N$ is called a} standard basis, \emph{ if $e_1,e_2 \in N_0$,
 and if the relations $Ve_i=\Pi e_i$ for $i=1,2$ and $e_3=Ve_1$ and $e_4=p^{-1}Ve_2$ hold.}
\end{Def}
It follows from the proof of Lemma \ref{critical} ({\em ii}) that there is a standard basis of $N$.

  \begin{Lem}\label{stb}
  Any superspecial $O_B$-lattice is the $W$-span of some standard basis of $N$.
  \end{Lem}
 {\em Proof.} Let $L$ be a superspecial $O_B$-lattice. Then $L_i=L_i^{\Pi^{-1}V}\otimes_{\Zp}W$ for $i=0,1$ (see the proof of Lemma \ref{critical} ({\em ii})). It follows from Lemma \ref{lemmain} that the inclusion $VL_0^{ \Pi^{-1}V} \subset L_1^{ \Pi^{-1}V}$ has index $1$. By the elementary divisor theorem we can choose a basis $e_3,e_4$ of $L_1^{ \Pi^{-1}V}$ such that $e_3, pe_4$ is a basis of $VL_0^{ \Pi^{-1}V}$. Setting $e_1=V^{-1}e_3$ and $e_2=V^{-1}pe_4$ we get the desired basis $e_1,e_2,e_3,e_4$. \qed
 \newline
 
  Denote by $\begin{rm}End\end{rm}(N,V)$ the space of endomorphisms of
$N$ which commute with $V$. Following \cite{KR1} we call an element $j \in \begin{rm}End\end{rm}(N,V)$ \emph{special} if it commutes with the $O_B$-action and its  trace is $0$. We denote the space of special endomorphisms by $V$.
Note that, by restricting an element $y\in V$ to the fixed ($\mathbb{Q}_p$-)module 
$N_0^{V^{-1}\Pi}$, we can identify $V$ with $M_2(\mathbb{Q}_p)^0$, the space of traceless matrices in $M_2(\mathbb{Q}_p)$, cf. \cite{KR1}, (2.2).

Let $*$ be a fixed $\mathbb{Q}_p$-linear automorphism of order $2$ of $B$. By the theorem of Skolem-Noether there is some $b_* \in B$ such that $*=\Int(b_*)$, i.e. $x^*=b_*xb_*^{-1}$. The element $b_*$ is unique up to multiplication by some element of $\mathbb{Q}_p^{\times}$. 
\begin{Def}
\emph{
An element $j\in\begin{rm}End\end{rm}(N,V)$ is \emph{$*$-special} if the following conditions are satisfied.
\begin{enumerate}[1.)]
\item
 $j\iota(a)j^{-1}=\iota(a^*) \ \forall a \in B$, 
 \item
  $s_j:=\iota(b_*^{-1})j$ has trace $0$.
\end{enumerate}
}
\end{Def}

Note that this definition is equivalent to the condition that $s_j$ is special since the first condition is equivalent to the condition that $s_j$ commutes with the $O_B$-action. Note also that the second condition is independent of the choice of $b_*$ since $b_*$ is unique up to multiplication by an element of $\mathbb{Q}_p^{\times}$.    We denote the space of $*$-special endomorphisms by $V[*]$. 

Since $*$ has order $2$, it follows that $b_*^2 \in \mathbb{Q}_p$. We will always choose $b_*$ in such a way that $b_*^2$ has valuation $0$ or $1$. Precisely one of these possibilities can be realized. In the first case we call $*$ unramified (since $\mathbb{Q}_p(b_*)$ is an unramified quadratic extension of $\mathbb{Q}_p$ in this case), in the second case we call $*$ ramified (since $\mathbb{Q}_p(b_*)$ is a ramified quadratic extension of $\mathbb{Q}_p$ in this case).
Hence $b_*$ is a unit in $O_B$ if $*$ is unramified, and $b_*=\Pi\cdot \varepsilon_*$, where $\varepsilon_* $ is a unit in $O_B$, if $*$ is ramified.  In the latter case we write $b_*^2=\eta_* p$ for some $\eta_*\in \Zp^{\times}$.

%-----------------------------------------------------------------------------------------------------------------------
%-----------------------------------------------------------------------------------------------------------------------

\section{Antispecial cycles}

Before we come to the definition of antispecial cycles we recall some definitions and statements of 
\cite{KR1}, §1.

A \emph{special formal $O_B$-module } over a $W$-scheme $S$ is a $p$-divisible formal group $X$ over $S$ 
of dimension 2 and height 4, with an $O_B$-action $\iota: O_B \rightarrow \End_S(X)$ such that the induced 
$\mathbb{Z}_{p^2}\otimes \mathcal{O}_S$-module $\Lie X$ is, locally on $S$, free of rank 1. 
We fix a special formal $O_B$-module $\mathbb{X}$ over $\Spec k$. Let us consider the following functor 
$\mathcal{M}$ on the category $\begin{rm} Nilp \end{rm}$ of $W$-schemes $S$ such that $p$ is locally 
nilpotent
in $\mathcal{O}_S$. It associates to a scheme $S \in \Nilp$  the set of isomorphism classes of  pairs
 $(X, \varrho)$
consisting of a special formal $O_B$-module $X$ over $S$ and an $O_B$-linear quasi-isogeny of height zero,
\[
\varrho: \mathbb{X}\times_{\Spec k}\overline{S} \rightarrow X\times_S \overline{S},
\]
where $\overline{S}=S\times_{\Spec W} \Spec k$. The functor $\mathcal{M}$ is representable by the
Deligne-Drinfeld formal 
scheme
$\widehat{\Omega}\times_{\begin{rm}Spf \end{rm}\ \mathbb{Z}_p}\begin{rm}Spf\end{rm}\ W$. Denote by
$\mathcal{B}=\mathcal{B}(PGL_2(\mathbb{Q}_p))$
the Bruhat-Tits building of $PGL_2(\mathbb{Q}_p)$. The formal scheme $\widehat{\Omega}$ is obtained by 
glueing formal open subschemes $\widehat{\Omega}_{\Delta}$, where $\Delta$ runs over the simplices of 
$\mathcal{B}$. We 
will
only need to know $\widehat{\Omega}_{\Delta}$ for $\Delta =$ standard vertex and $\Delta =$ standard edge.
So let $\Delta=[\Lambda_0]$ be the homothety class of the standard lattice
\[
\Lambda_0=[e_1,e_2],
\]
where $[e_1,e_2]$ here denotes the $\mathbb{Z}_p$-span of the standard basis in $\mathbb{Q}_p^2$. Then 
\begin{equation}\label{gewgl}
\widehat{\Omega}_{\Lambda_0}=\Spf \ \mathbb{Z}_p[T,(T^p-T)^{-1}]^{\wedge}.
\end{equation}
Here, $^{\wedge}$ denotes the $p$-adic completion.
If $\Delta =\Delta_0 =([\Lambda_0],[\Lambda_1])$ is the standard edge corresponding to $\Lambda_0=[e_1,e_2]$ and
$\Lambda_1=[pe_1,e_2]$, then 
\begin{equation}\label{ssgl}
\widehat{\Omega}_{\Delta_0}=\Spf \ \Zp[T_0,T_1,(1-T_0^{p-1})^{-1}, (1-T_1^{p-1})^{-1}]^{\wedge}/(T_0T_1-p).
\end{equation}

Any $k$-valued point of $\mathcal{M}$ corresponds to a special Dieudonn\a'e module with an $O_B$-action 
(as defined in section $1$).
We may choose $\mathbb{X}$ in its isogeny class so that its Dieudonn\a'e module $\mathbb{L}$ is
 $O_B$-superspecial. By Lemma \ref{stb} we therefore  find a standard basis 
$e_1,e_2, e_3=\Pi e_1, e_4=p^{-1}\Pi e_2$
 of $\mathbb{L}$. We suppose that 
 the isocrystal $N$ considered in 
 the preceding section equals 
the isocrystal of $\mathbb{L}$. Then any $k$-valued point of $\mathcal{M}$ corresponds to a $O_B$-lattice (in
$N$)
defined as
 above. The superspecial points (i.e. those $k$-valued points of $\mathcal{M}$ whose $O_B$-lattice is
 $O_B$-superspecial) are in one-to-one
 correspondence with the
 edges in $\mathcal{B}$.
  This correspondence can be chosen in such a way that $\mathbb{L}$ corresponds to the
 standard edge $\Delta_0$ defined above. In the formal scheme
$\widehat{\Omega}\times_{\begin{rm}Spf\end{rm}\ \mathbb{Z}_p}\begin{rm}Spf\end{rm}\ W$ the point 
$pt_{\Delta_0}$ 
corresponding to
$\Delta_0$ lies in $\widehat{\Omega}_{\Delta_0}\times_{\begin{rm}Spf\end{rm}\ \mathbb{Z}_p}\begin{rm}Spf\end{rm}\ W$ and is given there by the equations
 $T_0=T_1=0$. Any ordinary $k$-valued point of
$\mathcal{M}$ (i.e. whose $O_B$-lattice is  $O_B$-ordinary)
corresponds in $\widehat{\Omega}\times_{\begin{rm}Spf\end{rm}\ \mathbb{Z}_p}\begin{rm}Spf\end{rm}\ W$ to a
$k$-valued point of some $\widehat{\Omega}_{\Lambda}\times_{\begin{rm}Spf\end{rm}\ \mathbb{Z}_p}\begin{rm}Spf\end{rm} \ W$ for some vertex  $[\Lambda]\in \mathcal{B}$.

Since the isocrystal of $\mathbb{X}$ equals $N$ we can make the following
\begin{Def}
\emph{Let $j \in \End(N;V)$. Then the  \emph{ cycle $Z(j)$
associated to $j$} is the closed formal subscheme of $\mathcal{M}$ consisting of all points $(X,\varrho)$ such
that $\varrho \circ j \circ \varrho^{-1}$ lifts to an endomorphism of $X$.
If $j$ is a special endomorphism,  then  $Z(j)$ 
is called a \emph{special cycle} cf. \cite{KR1}, Definition 2.1. Let $*$ be a $\mathbb{Q}_p$-linear automorphism of order $2$ of $B$. A \emph{$*$-special cycle} is a cycle of the form  $Z(j)$ for some $j$ which is $*$-special. An \emph{antispecial cycle} is a cycle of the form $Z(\iota(\Pi)y)$ for some special endomorphism $y$.}
\end{Def}

The fact that $Z(j)$ is a closed formal subscheme of $\mathcal{M}$ follows from \cite{RZ}, Proposition 2.9, see also  \cite{KR1}, p. 167.
\begin{Rem}\label{maleps} 
\emph{Let $\varepsilon \in O_B^{\times}$, and let $j, \tilde{j}\in \End(N;V)$ be such that  
$j=\iota(\varepsilon) \tilde{j}$. Then  since 
 $\iota(\varepsilon)$ is invertible, it follows that $Z(\tilde{j})=Z(j)$.}
\end{Rem}

Suppose we are given  a $\mathbb{Q}_p$-linear automorphism $*$ of order $2$ of $B$. If $*$ is unramified  and $j$  is $*$-special, then it follows from Remark \ref{maleps} that $Z(j)=Z(s_j)$, hence $Z(j)$ is a special cycle. (Recall that $s_j=\iota(b_*^{-1})j$ and $b_*$ is a unit in $O_B$.) Conversely, it follows again from Remark \ref{maleps} that every special cycle is $*$-special for arbitrary unramified $*$. Since equations for  special cycles are known from  \cite{KR1}, in the sequel we will exclude the case that $*$ is unramified. Therefore, in the sequel we assume that $*$ is ramified.  Recall that $b_*=\Pi\cdot \varepsilon_*$ where $\varepsilon_*\in O_B^{\times}$ in this case. It follows again from Remark \ref{maleps} that for arbitrary ramified $*$ a cycle is $*$-special if and only if it is of the form $Z(\iota(\Pi) y)$ for some special endomorphism $y$, i.e. it is an antispecial cycle.  In particular it follows that the notion of $*$-special cycles does not depend on the particular choice of the (ramified) $*$. We now fix a ramified  $*$ for this section.

 Let $j \in V[*]$. By Remark \ref{maleps} we have $Z(j)=Z(\iota(\Pi) s_j)$.  A $k$-valued point of $\mathcal{M}$ belongs
 to $Z(j)$ if and
only if the $O_B$-lattice corresponding to that point is mapped by $j$ (or, equivalently, by $\iota(\Pi) s_j$) into itself. If $L$ is an 
$O_B$-lattice, this means $\iota(\Pi) s_j L \subset L$ and implies 
 $ps_j^2L\subset L$. From this we see $\nu_p(\begin{rm}  det \end{rm}(s_j)) \geq -1$, where $\nu_p$ is the
 valuation associated to $p$. So, in
 investigating $Z(j)$ we may assume $\nu_p(\begin{rm}  det \end{rm}(s_j)) \geq -1$, otherwise
 $Z(j)=\emptyset$.

Now, given a $*$-special endomorphism $j$ and an $O_B$-lattice $L=L_0 \oplus L_1$, 
we want to investigate
under which conditions the inclusion $jL\subset L$ holds.
Since $\Pi N_0=N_1$ and $\Pi N_1=N_0$, we see that $jL \subset L$ holds if and only if
 $jL_0\subset L_1$ and $jL_1\subset L_0$. 

 If $0$ is $O_B$-critical we set $A_0=L_0$ and $A_1=VL_1\subset A_0$. These are lattices in $N_0$.
The condition $jL\subset L$ then translates into the conditions $V\iota(\Pi)  s_j L_0 \subset VL_1$ and 
$\iota( \Pi ) s_j VL_1 \subset VL_0$, hence 

\begin{equation}\label{bedd1}
jL\subset L \ \Longleftrightarrow \ ps_jA_0 \subset A_1  \text{ and } 
s_j A_1 \subset A_0.
\end{equation}

If $1$ is $O_B$-critical we set $A_0=VL_0$ und $A_1=L_1$. We then get an analogous condition to (\ref{bedd1}) 
where the roles of 0 and 1 are interchanged.

 Now we assume, for example, that 0 is $O_B$-critical.
 By the elementary divisor theorem, there is a $W$-basis $f_1, f_2$ of $A_0$ for which $pf_1, f_2$ is a basis
 of $A_1$.
 
  If $ps_jA_0\subset A_1$, we can write
 \begin{equation}\label{bedd2}
 ps_j f_1=paf_1+cf_2 \text{ where } a,c \in W, 
 \end{equation}
and if $s_j A_1 \subset A_0$, we can write
\begin{equation}\label{bedd3}
 s_j f_2=bf_1+df_2 \text{ where } b,d \in W. 
 \end{equation}
 
 Conversely, if (\ref{bedd2}) and (\ref{bedd3}) are satisfied, so is the condition of (\ref{bedd1}).
  So, given $L$ 
 (where $0$ is $O_B$-critical) and 
 hence given $A_0$ and $A_1$,
 we have $jL\subset L$ if and only if there is a $W$-basis $f_1, f_2$ of $A_0$ such that (\ref{bedd2}) and 
 (\ref{bedd3}) hold.
The case where $1$ is $O_B$-critical is treated in the same way. 
(One just has to replace $A_0$ by $A_1$  and conversely in the above reasoning.)

 Thus we have shown the following
\begin{Lem}\label{inklemma}
Using the notations just introduced and assuming that $i \in \mathbb{Z}/2$ is $O_B$-critial for $L$ we have 
$jL\subset L$ if and only if there is a $W$-basis $f_0,f_1$ of $A_i$ for which the matrix expressing $s_j$ in
$f_0,f_1$ has the form
\begin{align}\label{matform}
s_j=
\begin{pmatrix}
 a& b \\ 
 p^{-1}c &d
\end{pmatrix}
\text{ where } a,b,c,d \in W.
\end{align}
\end{Lem}
\qed
\begin{Pro} \label{kwertig}
Let $j$ be an $*$-special endomorphism, where 
 $\nu_p(\begin{rm}  det \end{rm}(s_j)) \geq -1$.
Then, regarding $ps_j$ as a special endomorphism, we have an equality of $k$-valued points
\[
Z(j)(k)=Z(ps_j)(k).
\]
\end{Pro}

\emph{Proof.} It is clear that $Z(j)=Z(\Pi s_j)\subset Z(ps_j)$, hence $Z(j)(k)\subset Z(ps_j)(k)$.

Conversely, given a $k$-valued point $(X,\varrho)$ of $Z(ps_j)$,
 let $L=L_0\oplus L_1$ be its Dieudonn\a'e module which we can regard as an $O_B$-lattice in $N$.
 Assume (for example) that 0 is $O_B$-critical.
  Since $ps_j L \subset L$, we  have
 inclusions 

\begin{equation}\label{inkl}
ps_jL_0 \subset L_0 \text{ and }
ps_jVL_1 \subset VL_1.
\end{equation}

Since the inclusion $VL_1\subset L_0$ has index 1, there is a $W$-basis $f_1, f_2$ of $L_0$ for which 
$pf_1, f_2$
is a $W$-basis of $VL_1$. Because of (\ref{inkl}) we see that in the basis $f_1,f_2$ the
 endomorphism $ps_j$ has matrix 
\begin{align*}
ps_j=\begin{pmatrix}
 \tilde{a}& pb \\ 
 c &-\tilde{a}
\end{pmatrix}
\text{ where } \tilde{a},b,c \in W,
\end{align*}

and hence $\begin{rm}  det \end{rm}(s_j)=-p^{-2}\tilde{a}^2-p^{-1}bc$. But
 $\nu_p(\begin{rm}  det \end{rm}(s_j)) \geq -1$, so we conclude that $\tilde{a}$ is divisible by
 $p$. Therefore the matrix of $s_j$ can be written in the form 
\begin{align*}
s_j=\begin{pmatrix}
 a& b \\ 
 p^{-1}c &-a
\end{pmatrix}
\text{ where } a,b,c \in W.
\end{align*} From Lemma \ref{inklemma} we therefore see that $(X,\varrho)$  also belongs to $Z(j)$.
\qed
\newline

Given an $*$-special endomorphism $j$, we want to determine local equations for $Z(j)$. This will be done with
the help of the Grothendieck-Messing theory. 
We summarize below the facts from this theory which we will need.

Let $A$ be a local Artin ring with residue field $k$ (which is algebraically closed), and assume that $A$ is a
 $W$-algebra. Then
$p\in A$ is nilpotent. Assume that the  maximal ideal $I$ of $A$  (which is nilpotent)  carries a
nilpotent $pd$-structure (in the sense of \cite{M}, chapter III, Definition 1.1).
 Let $X_0$ be a $p$-divisible group over $k$ of height 4 and dimension 2.
Denote by $L$ the Dieudonn\a'e module of $X_0$, and define $P(X_0)_k=L\otimes_W k$ and 
$P(X_0)_A = L \otimes_W A$.
Let $\mathcal{F}_k=VL/pL \subset P(X_0)_k$. This gives us the
Hodge filtration
\begin{align}\label{hok}
0 \rightarrow \mathcal{F}_k \rightarrow P(X_0)_k \rightarrow L/VL \rightarrow 0.
\end{align}
To a lifting $X$ over $A$ of $X_0$ corresponds a lifting of the Hodge filtration (of $A$-modules),
\begin{align}\label{hoA}
0 \rightarrow \mathcal{F} \rightarrow P(X_0)_A \rightarrow \begin{rm}Lie\end{rm}(X) \rightarrow 0,
\end{align}
where $\mathcal{F}$ is a direct summand of $P(X_0)_A$ of rank 2, and where (\ref{hoA}) lifts (\ref{hok}).

This establishes an equivalence of categories between the category of liftings of $X_0$ over $k$ 
to some $X$ over $A$ and the category of
liftings of  the Hodge filtration (\ref{hok}) to filtrations of the form (\ref{hoA}).

An endomorphism $\phi : X_0 \rightarrow X_0$  gives rise to  an endomorphism
 $\Phi : L \rightarrow L$. Then
 $\phi$ lifts to an endomorphism of $X$ if and only if
the endomorphism induced by $\Phi$ on $P(X_0)_A$ maps the submodule $\mathcal{F}$ into itself. 
(In this situation
we will denote the induced endomorphism by $\Phi$ as well.)
Now assume that $X_0$ is equipped with an $O_B$-action $\iota:O_B \rightarrow \End_k(X_0)$. Then we get an
$O_B$-action  $\iota:O_B \rightarrow \End(L)$. We apply the equivalence of categories just mentioned to get
the following 
\begin{Pro}\label{grome}
Using the same notations as above, the following categories are equivalent
\begin{enumerate}[(i)]
\item The category of liftings of $X_0$ 
to  $A$,  also lifting the $O_B$-action
\item The category of
liftings of  the Hodge filtration (\ref{hok}) to $P(X_0)_A$ 
 which are stable under
the induced $O_B$-action. 
\end{enumerate}
\qed
\end{Pro}
 
\begin{Lem}\label{pdringe}
The rings $A=W/p^n$ for $n \in \mathbb{N}$ and if $p>3$ also the rings $A=W[x]/(x^2-p \varepsilon, x^r)$, 
where $\varepsilon \in W^{\times}$ and $r\geq
2$ satisfy the requirements of the Grothendieck-Messing theory mentioned above, i.e., they are local Artin
rings and $W$-algebras  with residue field $k$,  and  their maximal ideals carry 
a nilpotent $pd$-structure.
\end{Lem} 
\emph{Proof.}
i) $A=W/p^n$ is obviously a local Artin ring and a $W$-algebra in the natural way. That $I=(p)$ carries a
nilpotent $pd$-structure can be seen as follows: 
If we denote by $\nu_p$ the valuation associated to $p$ we have for any positive integer $i$,
\[
\nu_p(i!)=\sum\nolimits_{k=1}^{\infty}[\frac{i}{p^k}] <\sum\nolimits_{k=1}^{\infty}\frac{i}{3^k} =
 \frac{i}{2} < i= \nu_p(p^i).
\]
(Here, [ ], denotes Gauss-brackets. The first inequality uses $p \geq 3$.) This shows that for any $z\in I$ 
we have a well-defined expression $\frac{z^k}{k!}$ and that there is some $N$ such that, whenever
$n_1+...+n_l\geq N$ we have $\frac{z^{n_1}}{n_1!}\cdot ... \cdot \frac{z^{n_l}}{n_l!}=0$ for all $z\in I$.
 Hence we see that $I$
carries a nilpotent $pd$-structure.

ii) $A=W[x]/(x^2-p \varepsilon, x^r)$ clearly is a local Artin ring and a $W$-algebra in the canonical way.
 Let $\nu_x$ denote the valuation associated to $x$. To see that $I=(x)$ carries a nilpotent $pd$-structure 
we consider $\nu_x(i!)$ for any
positive integer $i$. Now, $\nu_x(i!)=2 \cdot \nu_p(i!)$, and since $p \geq 5$, 
\[
\nu_x(i!)=2\cdot\nu_p(i!)=2\cdot \sum\nolimits_{k=1}^{\infty}[\frac{i}{p^k}] 
<2\cdot \sum\nolimits_{k=1}^{\infty}\frac{i}{5^k} 
= \frac{i}{2}.
\]
Hence we see as above that there is a nilpotent $pd$-structure on $I$.
\qed

\begin{Pro}\label{wmodpn}
Let $n \in \mathbb{N}$, and let $j$ be a $*$-special endomorphism with
$\nu_p(\begin{rm}  det \end{rm}(s_j)) \geq -1$. Then $Z(j)$ and $Z(ps_j)$ 
have the same $W/p^n$-valued points, 
\[
Z(j)(W/p^n)=Z(ps_j)(W/p^n).
\]
\end{Pro}
\emph{Proof.} The case $n=1$ is treated in Proposition \ref{kwertig}, so we may assume that $n \geq 2$. 
Since $Z(j)=Z(\Pi s_j)\subset Z(p s_j)$, it is enough to show that $Z(p s_j)(W/p^n)\subset Z(\Pi s_j)(W/p^n)$.
 We fix a $k$-valued point $(X_k,\varrho_k)$ of
$Z(ps_j)$ and consider liftings to $W/p^n$ with the Grothendieck-Messing theory. Let $L$ be the Dieudonn\a'e
 module of  $(X_k,\varrho_k)$. Let $P=L/p^nL$. Then the Hodge filtration of $(X_k,\varrho_k)$ is given by
 $VL/pL\subset L/pL$. 
 By Proposition \ref{grome} and Lemma \ref{pdringe} to a
  lifting $(X,\varrho)$ of  
 $(X_k,\varrho_k)$ over $W/p^n$ (within $Z(ps_j)$) corresponds a lifting of 
the Hodge filtration $\mathcal{F}\subset P$ stable under the $O_B$-action and under $ps_j$.

 From the action of $\mathbb{Z}_{p^2}$ we obtain  $\mathbb{Z}/2$-gradings,
\[
\mathcal{F}=\mathcal{F}_0 \oplus \mathcal{F}_1
\]
and 
\[
P=P_0\oplus P_1.
\]

Let $\mathcal{F}_0=<\overline{f_0}>$ and $\mathcal{F}_1=<\overline{f_1}>$ for suitable elements $\overline{f_i} \in P_i$. 

\begin{bf} Claim:\end{bf} \emph{ For some units
$\varepsilon_0, \varepsilon_1$ we have either 
\begin{equation}\label{pmal}
\Pi\overline{f_0}= \varepsilon_1 p \overline{f_1} \text{ and } \Pi\overline{f_1}= \varepsilon_0  
\overline{f_0},
\end{equation}
or
\begin{equation}\label{1mal}
\Pi\overline{f_0}= \varepsilon_1  \overline{f_1} \text{ and } \Pi \overline{f_1}= \varepsilon_0 p 
\overline{f_0}.
\end{equation}}
To see this, choose preimages $f_0 \in L_0$ of $\overline{f_0}$ and $f_1\in L_1$ of $\overline{f_1}$. 
Since $\Pi L_i \subset L_{i+1}$, we have
\[
\Pi f_0 = \varepsilon_1 p^{\nu_1}f_1 +p^n\xi_1 \text{ and } \Pi f_1 = \varepsilon_0 p^{\nu_0}f_0+p^n\xi_0,
\]
for some integers $\nu_0,\nu_1$, some units $\varepsilon_0, \varepsilon_1$ and some $\xi_0,\xi_1 \in L$.
Therefore, for some $\Xi \in L$, we have 
\[
pf_0=\varepsilon_0\varepsilon_1p^{\nu_0+\nu_1}f_0+p^n\Xi
\]
and hence
(using $n \geq 2$) we get $\nu_0+\nu_1=1$, i.e. either $\nu_0=0$ and $\nu_1=1$ or $\nu_0=1$ and $\nu_1=0$.
This confirms the claim. 

Because of the symmetry of (\ref{pmal}) and (\ref{1mal}), 
we may, for example, assume that (\ref{pmal}) holds.
 By replacing $f_0$ by $\Pi f_1$ and thereby
changing $\overline{f_0}$ only by a unit and not changing $\mathcal{F}$ we can  assume that
\[
\Pi f_0 = pf_1 \text{ and } \Pi f_1 = f_0.
\]
We choose $\overline{l_0}\in P_0$ with preimage $l_0\in L_0$ and $\overline{l_1}\in P_1$ with preimage 
$l_1\in L_1$ such
that $P=<\overline{l_0},\overline{f_0},\overline{l_1},\overline{f_1}>$. By Nakayama's lemma we also have
 $L=<l_0,f_0,l_1,f_1>$.
Let 
\[
\Pi l_0 = al_1 + bf_1\text{ } \text{ and } \text{ } \Pi l_1 = cl_0+df_0\text{ } \text{ where } \text{ }
a,b,c,d \in W,
\]
then
\[
pl_0=a(cl_0+df_0)+bf_0,
\]
hence $ac=p$ and $ad+b=0$. 

\begin{bf} Claim:\end{bf} \emph{ $a$ is a unit.}

 To see this, assume that $a$ is divisible by $p$. Then  $b$ is also divisible by $p$.
Therefore, 
\[
\Pi L_0 \ = \ <\Pi f_0 , \Pi l_0>\ = \ <pf_1,  p(\frac{a}{p}l_1+ \frac{b}{p}f_1)>.
\]
But this contradicts the fact that $L_1/ \Pi L_0$ is a $k$-vector space of dimension $1$. Therefore $a$ is a 
unit which proves the claim.

By Proposition \ref{kwertig} we have $\Pi  s_j L \subset L$. Writing 
\[
s_j f_0=rl_0+sf_0
\] we may assume $\nu_p(r), \nu_p(s)\geq -1$, where $\nu_p$ is the valuation associated to $p$. (This follows
>>from $\Pi s_j L \subset L$ and hence $ps_jL\subset L$.) We now want to show that $\Pi  s_j \mathcal{F} \subset
\mathcal{F}$. If $l \in L$, denote by $\overline{l}$ its image in $P$. First,
\[
\Pi  s_j \overline{f_0}=\overline{\Pi  s_j f_0}=\overline{s_j  \Pi f_0}=\overline{s_jp f_1}
=ps_j\overline{f_1}\in \mathcal{F}, 
\]
since $ps_j\mathcal{F}\subset \mathcal{F}$. Therefore,
\[
\mathcal{F} \ni \Pi  s_j\overline{f_0}=\overline{\Pi  s_j f_0}=\overline{\Pi(rl_0+sf_0)}=
\overline{r(al_1+bf_1)+psf_1}.
\]
Now, $r(al_1+bf_1)+psf_1\in L$, and $a$ is a unit, hence $\nu_p(r)\geq 0$ and further, since 
$\overline{r(al_1+bf_1)+psf_1}\in \mathcal{F}$, and $a$ is a unit, we even have 
$r \equiv 0 \ \begin{rm} mod \end{rm} \ p^n$.
Now, $L \ni \Pi s_j f_1 = s_jf_0$, so $\nu_p(s)\geq 0$. Since $r \equiv 0 \ \begin{rm} mod \end{rm}\ p^n$,
we get 
\[
\overline{\Pi  s_j f_1} =\overline{s_j \Pi f_1} =\overline{ s_jf_0}=\overline{sf_0}\in \mathcal{F}.
\]
This completes the proof.
\qed
\newline

In the following statement we denote by an upper index "ord" the intersection with the ordinary locus of
$\mathcal{M}$ resp. $\widehat{\Omega}$, i.e. the open formal subscheme formed by the complement of the
superspecial points.

\begin{Pro}\label{ordloc}
Let $j$ be a $*$-special endomorphism with $j^2 \neq 0$, where 
 $\nu_p(\begin{rm}  det \end{rm}(s_j)) \geq -1$.
Then 
\[
Z(j)^{\begin{rm}ord\end{rm}}=Z(ps_j)^{\begin{rm}ord\end{rm}}.
\]
\end{Pro}
\emph{Proof.} It is clear that $Z(j)^{\begin{rm}ord\end{rm}}\subset Z(ps_j)^{\begin{rm}ord\end{rm}}$.
Let $[\Lambda]$ be a lattice in $\mathbb{Q}_p^2$. In order to show
$Z(j)\cap (\widehat{\Omega}_{[\Lambda]}\times_{\begin{rm}Spf\end{rm}\ \mathbb{Z}_p}\begin{rm}Spf\end{rm}\ W)
=Z(ps_j)\cap (\widehat{\Omega}_{[\Lambda]}\times_{\begin{rm}Spf\end{rm}\ \mathbb{Z}_p}\begin{rm}Spf\end{rm}\ W)$ we
may assume that $[\Lambda]=[\Lambda_0]$ is the standard lattice, cf. \cite{KR1}, Proposition 3.2. We write in the 
basis $e_1,e_2$ of $N_0$
\[
ps_j=
\begin{pmatrix}
 a& b \\ 
c &-a
\end{pmatrix}
= 
p^m \cdot
\begin{pmatrix}
 \overline{a}& \overline{b} \\ 
 \overline{c} &-\overline{a}
\end{pmatrix},
\]
where $\overline{a},\overline{b},\overline{c}\in \mathbb{Z}_p$ are not simultaneously divisible by $p$. The equation for
$Z(ps_j)\cap (\widehat{\Omega}_{[\Lambda]}\times_{\begin{rm}Spf\end{rm}\ \mathbb{Z}_p}\begin{rm}Spf\end{rm}\ W)$ can
be written as 
\[
p^m \cdot (\overline{b}T^2-2\overline{a}T-\overline{c})=0,
\]
see \cite{KR1}, Proposition 3.2. 
Let $Z(j)\cap (\widehat{\Omega}_{[\Lambda]}\times_{\begin{rm}Spf\end{rm}\ \mathbb{Z}_p}\begin{rm}Spf\end{rm}
\ W)$
be described by the Ideal $I$  of $W[T,(T^p-T)^{-1}]^{\wedge}$ as a closed subscheme of 
$\widehat{\Omega}_{[\Lambda]}\times_{\begin{rm}Spf\end{rm}\ \mathbb{Z}_p}\begin{rm}Spf\end{rm}\ W=
\begin{rm}Spf\end{rm}\ W[T,(T^p-T)^{-1}]^{\wedge}$ (comp. (\ref{gewgl})).

 \begin{bf} Claim:\end{bf} \emph{ Every element of $I$ is divisible by $(\overline{b}T^2-2\overline{a}T-\overline{c})$.}

If 
\[
\begin{rm}Spf\end{rm}\ W[T,(T^p-T)^{-1}]^{\wedge}/(\overline{b}T^2-2\overline{a}T-\overline{c})
\]
is not empty, then 
 (by \cite{KR1}, Proposition 3.2)  
it meets the special fibre of
$Z(ps_j)$ in two different ordinary points. 
 From the form of the matrix in (\ref{matform})  (for some basis of $\Lambda$), we see that the rank of the 
 matrix obtained 
 from $ps_j$ by reduction 
mod $p$ is at most 1. On the other hand, in \cite{KR1}, Proposition 2.3 it is shown that if $Z(ps_j)(k)\cap
(\widehat{\Omega}_{[\Lambda]}\times_{\begin{rm}Spf\end{rm}\ \mathbb{Z}_p}\begin{rm}Spf\end{rm}\ W) \neq
\emptyset$, 
then this rank equals $0$ or $2$. Hence we conclude that $ps_j$ is divisible by $p$.
Hence $s_j(\Lambda)\subset \Lambda$. Therefore,  
\[
\begin{rm}Spf\end{rm}\ W[T,(T^p-T)^{-1}]^{\wedge}/(\overline{b}T^2-2\overline{a}T-\overline{c})
\subset Z(s_j)\subset Z(j),
\]
as asserted.

\begin{bf} Claim:\end{bf} \emph{
Every element of $I$ is divisible by $p^m \cdot(\overline{b}T^2-2\overline{a}T-\overline{c})$.}

Suppose there is an element $Q\in I$ which is not divisible by
$p^m \cdot(\overline{b}T^2-2\overline{a}T-\overline{c})$. By multiplying with a suitable power of $p$ we 
may suppose that $Q$ is divisible by
$p^{m-1}(\overline{b}T^2-2\overline{a}T-\overline{c})$ and write 
\[
Q=q\cdot p^{m-1} \cdot(\overline{b}T^2-2\overline{a}T-\overline{c}),
\]
 where we may suppose that $q$ is a polynomial in $T$
in which no coefficient is divisible by $p$. (This assumption is allowed, since 
\[
(p^m \cdot (\overline{b}T^2-2\overline{a}T-\overline{c})) \subset I,
\]
which follows from the inclusion $Z(j)\subset Z(ps_j)$.)

 Let  $\tilde{q}$ be the image of
\[
 q\cdot (\overline{b}T^2-2\overline{a}T-\overline{c})
\cdot (T^p-T)
\]
 in $k[T]$. By Gauss' lemma $\tilde{q} \neq 0$, and we can find $\tau \in k^{\times}$ with
$\tilde{q}(\tau) \neq 0$. Choose a preimage $t \in W/p^m$ of $\tau$. Then 
$q(t)\cdot (\overline{b}t^2-2\overline{a}t-\overline{c})\cdot (t^p-t)$ is a unit in $W/p^m$.

Now we are ready to apply Proposition \ref{wmodpn}, which says that any $W/p^m$-valued point 
\[
\phi:W[T,(T^p-T)^{-1}]^{\wedge}/(p^m
\cdot(\overline{b}T^2-2\overline{a}T-\overline{c}))  \longrightarrow W/p^m
\]
  factors through 
$W[T,(T^p-T)^{-1}]^{\wedge}/I$. 

Define $\phi$ by $\phi(T)=t$ and to be $W$-linear. Then $\phi(Q)=Q(t) \neq 0$. But then $\phi$ does not
factor through $W[T,(T^p-T)^{-1}]^{\wedge}/I$. This contradiction proves the claim.

Altogether we get 
\[
(p^m \cdot(\overline{b}T^2-2\overline{a}T-\overline{c})) \ \subset \ I \  \subset \
( p^m \cdot(\overline{b}T^2-2\overline{a}T-\overline{c})),
\] 
  and therefore 
$Z(j)^{\begin{rm}ord\end{rm}}=Z(ps_j)^{\begin{rm}ord\end{rm}}$.
\qed
\newline

The next proposition gives local equations of $Z(j)$ for a $*$-special endomorphism $j$ in a neighborhood of a superspecial point. Any
superspecial point $x$ corresponds to a simplex $\Delta = ( [ \Lambda ],[ \tilde{\Lambda} ])$. We may 
suppose
that $\Delta = \Delta_0 = ( [ \Lambda_0 ],[ \Lambda_1 ])$ is the standard simplex
 (see \cite{KR1}, Proposition 3.3.) 
We then have $x\in Z(ps_j)$ and,
 with respect to the basis $e_1, e_1$ of $N_0$, we write again
\[
ps_j=
\begin{pmatrix}
 a& b \\ 
c &-a
\end{pmatrix}
= 
p^m \cdot
\begin{pmatrix}
 \overline{a}& \overline{b} \\ 
 \overline{c} &-\overline{a}
\end{pmatrix}.
\]

Let $g$ be a special endomorphism with $g^2 \neq 0$ and $\nu_p(\begin{rm}  det \end{rm}(g)) \geq 0$. In \cite{KR1}, Proposition 3.3  local equations for the special cycle $Z(g)$ are given, 
where the following 
 three cases are distinguished:
\begin{enumerate}[(i)]

\item    $[ \Lambda_0 ]$ is strictly closer than $[ \Lambda_1 ]$ to the fixed point set $\mathcal{B}^g$ in the
 Bruhat-Tits-building $\mathcal{B}$.
 
\item The fixed point set $\mathcal{B}^g$ is the midpoint of $\Delta$, i.e. $[ \Lambda_0 ]$ and $[ \Lambda_1 ]$ both have 
 distance $1/2$ to $\mathcal{B}^g$.
 
\item $\Delta$ lies in the fixed apartement $\mathcal{B}^g$, i.e. $[ \Lambda_0 ]$ and $[ \Lambda_1 ]$ both have 
 distance $0$ to $\mathcal{B}^g$. 
 \end{enumerate}
 In the following proposition we follow these cases for the special endomorphism
 $ps_j$ to give local equations for $Z(j)$ (in case $j^2 \neq 0$ with $\nu_p(\begin{rm}  det \end{rm}(s_j)) \geq -1$) in a neighborhood of the superspecial point $x=pt_{\Delta}$.
  Using these notations and the description for
 $\widehat{\Omega}_{\Delta_0}$ in (\ref{ssgl}) we state
\begin{Pro}\label{acycglss}
Let $p>3$. Let $x=pt_{\Delta_0}$ be the superspecial point as above.
\begin{enumerate}[(i)]
\item Suppose $[ \Lambda_0 ]$ is strictly closer than $[ \Lambda_1 ]$ to the fixed point set $\mathcal{B}^{s_j}$.
We
then have $m \geq 1$, and $Z(j)$ is locally around $x$ given by the equation
\[
T_0\cdot p^{m-1}=0.
\]
\item Suppose $\mathcal{B}^{s_j}$ is the midpoint of $\Delta_0$. Then $\overline{b}$ is divisible by $p$, and  
 $Z(j)$ is locally around $x$ given by the
equation
\[
p^m\cdot (b_0T_0-2\overline{a}-\overline{c}T_1)=0,
\]
where $\overline{b}=p\cdot b_0$.
\item Suppose $\Delta_0$ lies in the fixed apartement $\mathcal{B}^{s_j}$. Then $m\geq 1$, and  $Z(j)$ is locally 
around $x$ given by the equation
\[
p^m=0.
\]
\end{enumerate}
\end{Pro}
\emph{Proof} of \emph{i}). We have $m\geq 1$, and there is an (affine) neighborhood  of $x$ in which $Z(ps_j)$ 
is
 given by the equations,
\[
T_0^2p^{m-1}=p^m=0,
\]
see \cite{KR1}, Proposition 3.3. (All references to \cite{KR1} in this proof refer to Proposition 3.3.) 
Let this neighborhood be given by 
\[
\begin{rm}Spf\end{rm}\ W[T_0,T_1,\gamma^{-1}]^{\wedge}/(T_0T_1-p),
\]
for some 
\[
\gamma \in W[T_0,T_1]^{\wedge} \setminus (T_0, T_1, T_0T_1-p)
\]
which is divisible by
 $(1-T_0^{p-1})(1-T_1^{p-1})$.
Then $Z(j)$ is in this  neighborhood  given by an ideal $I$ of the ring
\[
R:=W[T_0,T_1,\gamma^{-1}]^{\wedge}/(T_0T_1-p,T_0^2p^{m-1},\ p^m)
\]
as a closed subscheme of
 $Z(ps_j)$.
 
 Suppose there is an element $Q$ of $I$ which is not divisible by $p^{m-1}$. By multiplying with a suitable
 power of $p$ we may suppose then 
 that $Q$ is divisible by $p^{m-2}$ and write
\[
Q\  = \ p^{m-2}\cdot(\sum\nolimits_{i=0}^{l_0} a_i T_0^i\ + \ \sum\nolimits_{j=1}^{l_1}b_j T_1^j)\
  + \ p^{m-1}\cdot \eta,
\]
where the $a_i$ and $b_j$ are either units in $W$ or zero, and where $\eta \in R$. Let $q=Q/p^{m-2}$.

\begin{bf} Claim:\end{bf} \emph{ $a_i=0 \ \forall i$.}

 To see this, 
suppose  there is an $a_i \neq 0$. Let $\tilde{q}$ be the image of $q\cdot \gamma$ in $k[T_0, T_1]$. Since 
$\gamma  \notin (T_0, T_1, T_0T_1-p) = (T_0, T_1, p)$, we find  $\tau \in k^{\times}$ such that 
$\tilde{q}(\tau, 0)\neq 0$. Let $t$ be a lifting of $\tau$ in $W/p^{m-1}$. We are going to apply Proposition
\ref{wmodpn} again which here says that any $W/p^{m-1}$-valued point
\[
\phi : R \longrightarrow W/p^{m-1}
\]
factors through $R/I$. Now define $\phi$ by $\phi(T_0)=t, \ \phi(T_1)=p\cdot
t^{-1}$ and to be $W$-linear. Then $\phi$ is well defined ($\phi(\gamma)$ is a unit) but 
$
\phi(Q)=Q(t,p\cdot t^{-1}) \neq 0
$, since the image of $Q(t,p\cdot t^{-1})/p^{m-2}$ in $k$ equals $\tilde{q}(\tau, 0)$. Therefore $\phi$ does not factor
 through
 $R/I$. This contradiction confirms the claim.
 
  In the same way one sees that $b_j=0 \ \forall j$. This shows that
 any element of $I$ is divisible by $p^{m-1}$. Now we show that any element of $I$ is divisible by
  $T_0\cdot p^{m-1}$. Any $Q\in I$ can be written as 
  \[
 Q\ = \  a\cdot p^{m-1}T_0 \ + \ p^{m-1} \cdot \sum\nolimits_{j=0}^{l}b_j T_1^j,
  \]
where $a$ and  all $b_j$ are either units in $W$ or zero. Let $q=Q/p^{m-1}$. 

\begin{bf} Claim:\end{bf} \emph{ $b_j=0 \ \forall j$}.

Assuming there is a $b_j \neq 0$, we find
again $\tau \in k^{\times}$ such that $(q\cdot \gamma )(0, \tau)$ does not vanish in $k$. We lift $\tau$ to 
$t\in W/p^m$ and define a $W$-linear $W/p^m$-valued point $\phi : R \rightarrow W/p^m$ by 
$\phi(T_0)=p\cdot t^{-1}$ and $\phi(T_1)=t$. We see as above that it does not factor through $R/I$. This
contradiction shows that  all
$b_j$ vanish. Hence either $I=(p^{m-1} \cdot T_0)$ or $I=0$. To confirm the claim of \emph{(i)} it
 is thus enough to 
show that $I \neq 0$.

In order to show that $I \neq 0$, we show the existence of an $W[z]/(z^2-p,\ p^m)$-valued point of $Z(ps_j)$ 
whose underlying
$k$-valued point is $x$, and which does not belong to $Z(j)$.  By the Grothendieck-Messing theory, this can be
done by constructing a lifting of the Hodge filtration of $x$ 
over $W[z]/(z^2-p,\ p^m)$   which is stable under $ps_j$, but no 
stable under $j$.

Let $L$ be the Dieudonn\a'e module of $x$. Then a
basis of $L$ is given by $e_1, e_2, e_3= \Pi e_1,e_4= p^{-1}\Pi e_2$, since $x$ corresponds to the 
standard simplex
$\Delta_0$, i.e $L=\mathbb{L}$. Let 
\[ 
L_z=L\otimes_W W[z]/(z^2-p),
\] 
 and let 
 \[ 
P=L\otimes_W W[z]/(z^2-p,p^m)=L_z/p^m L_z.
\] 
Define 
\[
\begin{split} 
&f_0 \ = \ e_2  + z(e_1+e_2)\ \in L_z, \\
&f_1 \ = \ e_3  + pe_4 +ze_4 \ \in L_z.
\end{split}
\]
Denote by $\overline{f_0}$ resp. $\overline{f_1}$ the image of $f_0$ resp. $f_1$ in $P$, and define
a filtration $\mathcal{F}\hookrightarrow P$ by $\mathcal{F} \ = \ < \overline{f_0}, \overline{f_1}>$. 
To see that it lifts the Hodge filtration of $x$ we note that the image of $\mathcal{F}$ in
$L\otimes_W k$ equals the span of the images of $e_2,e_3$, and this is the image of $\Pi L= VL$. Further we
have $\Pi f_0=zf_1$ and $\Pi f_1=zf_0$ which shows that $\mathcal{F}$ is $O_B$-stable. 

The map $ps_j$ induces on $P$
the zero map, in particular $\mathcal{F}$ is stable under it. 
Now let us show that $\mathcal{F}$ is not stable under $\Pi s_j$.
For $l\in L_z$ let us denote its image in $P$ by
$\overline{l}$. A short calculation shows
\[
\overline{\Pi s_j f_0}=\overline{(p^{m-1}\overline{b} + zp^{m-1}\overline{a}+zp^{m-1}\overline{b})e_3}.
\]
 Let us suppose this is contained in $\mathcal{F}$.
Then we can find $r\in W[z]/(z^2-p)$ with $\overline{rf_1}=\overline{\Pi s_j f_0}$, hence 
$\overline{re_3}=\overline{\Pi s_j f_0}$ and $\overline{r(p+z)e_4}=0$
 and hence,
\[
\begin{split} 
r & \equiv  p^{m-1}\overline{b} + zp^{m-1}\overline{a}+zp^{m-1}\overline{b} \  \begin{rm}  mod \end{rm} \ p^m, \\
(p+z)r & \equiv 0  \ \begin{rm}  mod \end{rm} \ p^m. 
\end{split}
\]
But if we multiply the first congruence by $p+z$ , we see that the second is not fulfilled, since $\overline{b}$ 
is not divisible by $p$, see \cite{KR1}. Hence indeed
$\Pi s_j \mathcal{F} \nsubseteq \mathcal{F}$.
\newline

\emph{Proof} of \emph{ii}). We consider 
the (affine) neighborhood of $x$ given by  (\ref{ssgl})
\[
\begin{rm}Spf\end{rm}\ W[T_0,T_1,\gamma^{-1}]^{\wedge}/(T_0T_1-p),
\]
where
\[
\gamma = (1-T_0^{p-1})(1-T_1^{p-1}).
\]
By \cite{KR1} in this neighborhood $Z(ps_j)$ is described 
by the equations
\[
p^m T_0 (b_0T_0 -2\overline{a}-\overline{c}T_1)=p^m T_1 (b_0T_0 -2\overline{a}-\overline{c}T_1)=0.
\]
Then $Z(j)$ is given by an ideal $I$ of the ring
\[
R:=W[T_0,T_1,\gamma^{-1}]^{\wedge}/
(T_0T_1-p,\ p^m T_0 (b_0T_0 -2\overline{a}-\overline{c}T_1), \ p^m T_1 (b_0T_0 -2\overline{a}-\overline{c}T_1))
\]
as a closed subscheme of
 $Z(ps_j)$ in that neighborhood. 
 
 \begin{bf} Claim:\end{bf} \emph{Any element of $I$ is divisible by 
 $b_0T_0 -2\overline{a}-\overline{c}T_1$.}
 
  (In case $m\geq 1$  this is clear, since then $Z(s_j)\subset \ Z(j)$ and 
 the ideal of
 $Z(s_j)$ is contained in $(b_0T_0 -2\overline{a}-\overline{c}T_1)$.)  Let 
 \[
 S=W[T_0,T_1,\gamma^{-1}]^{\wedge}/ (T_0T_1-p, b_0T_0 -2\overline{a}-\overline{c}T_1).
 \]
Denoting the image of $I$ in $S$  by
$\overline{I}$, we must show that $\overline{I}=0$. From the relations $T_0T_1-p=0$ and  $b_0T_0 -2\overline{a}-\overline{c}T_1=0$
 and the fact that $b_0$ and $\overline{c}$ are units (see \cite{KR1}), we get $T_1= c^{-1} (b_0 T_0-2\overline{a})$ and 
 $T_0^2=2ab_0^{-1}T_0 + pc b_0^{-1}$, showing that any element in $S$ and so in particular 
 any $Q\in \overline{I}$ can be written
 in the form $Q=rT_0+s$, where $r,s\in W$. Now let $\varepsilon = (\frac{\overline{a}}{pb_0})^2p+\frac{c}{b_0}$ 
 (note that $\overline{a}$ is divisible by $p$, see \cite{KR1}) and define
  for any $n\in \mathbb{N}$ a $W$-linear homomorphism 
 \[
 \phi_n: S \longrightarrow W[z]/(z^2-p\varepsilon, p^n),\ \ \ T_0 \mapsto \frac{\overline{a}}{b_0}+z, \  \ \  T_1 
 \mapsto
 \frac{p}{\frac{\overline{a}}{b_0}+ z}.
 \]
 One easily checks that $\phi_n$ is well defined. Let $\pi:\  R \longrightarrow S$ be the canonical 
 projection.
 Consider for any $n$ the $W[z]/(z^2-p\varepsilon, p^n)$-valued point of $Z(ps_j)$,
\[
 \alpha_n= \phi_n \circ \pi: R \longrightarrow W[z]/(z^2-p\varepsilon, p^n).
\] 
These maps are compatible with varying $n$.
We want to show that $\alpha_n$ is also  a $W[z]/(z^2-p\varepsilon, p^n)$-valued point of $Z(\Pi s_j)$.
 Let $L=\mathbb{L}$ be the Dieudonn\a'e module of $x$, let 
$L_z= L\otimes_W  W[z]/(z^2-p\varepsilon)$, and let 
  $\mathcal{F}_n \hookrightarrow L_z/ (p^n L_z)$ be the Hodge filtration associated to $\alpha_n$. Then
 $\mathcal{F}_n=\mathcal{F}_{n+1}/(p^n L_z)$, since the $p$-divisible group corresponding to
 $\alpha_{n+1}$ lifts the one corresponding to $\alpha_n$. Let 
 $\mathcal{F}_{n+1}=<\overline{f_0},\overline{f_1}>$,
  where $\overline{f_i}$ is contained in
 the index $i$-part of $\mathcal{F}_{n+1}$. Further, for $i\in \mathbb{Z}/2$,
 let $f_i$ be a lifting of $\overline{f_i}$ in 
 $L_z$. Since $ps_j\mathcal{F}_{n+1}\subset \mathcal{F}_{n+1}$, we have
 $ps_jf_i=\lambda_i f_i + p^{n+1}\xi$ for some $\lambda_i \in W[z]/(z^2-p\varepsilon)$ and some
 $\xi \in L_z$. Therefore $p\Pi s_j f_i=\tilde{\lambda}f_{i+1}+p^{n+1}\tilde{\xi}$ for some 
 $\tilde{\lambda_i} \in W[z]/(z^2-p\varepsilon)$ and some 
 $\tilde{\xi} \in L_z$. Because of Proposition \ref{kwertig} we have $\Pi s_j L_z \subset L_z$, hence 
 $p\Pi s_j f_i\in pL_z$.
  Since $f_i$  belongs to a basis of $L_z$ (Nakayama's lemma),  we conclude that
 $\tilde{\lambda}$ is divisible by $p$, hence for some $\mu \in W[z]/(z^2-p\varepsilon)$ we have
 $\Pi s_j f_i= \mu f_{i+1}+p^n\tilde{\xi}$. Therefore $\Pi s_j \mathcal{F}_n \subset \mathcal{F}_n$, which shows 
 that $\alpha_n$ is a $  W[z]/(z^2-p\varepsilon, p^n)$-valued point of $Z(\Pi s_j)$.
 
 It follows that $\phi_n$ must factor through $S/\overline{I}$. Assume now that $\overline{I}\neq 0$. Then we find 
 $0 \neq rT_0+s \in\overline{I}$ and a sufficiently large $n$ such that $r(\frac{a}{b_0}+ z)+s$ does not
 vanish in $W[z]/(z^2-p\varepsilon, p^n)$. Therefore for such $n$
 the homomorphism $\phi_n$  will not factor through $S/\overline{I}$. This contradiction
  shows
 $\overline{I}=0$, hence any element of $I$ is divisible by $b_0T_0 -2\overline{a}-\overline{c}T_1$. 
 
 \begin{bf} Claim:\end{bf} \emph{ Any element of $I$ is divisible by $p^m(b_0T_0 -2\overline{a}-\overline{c}T_1)$}.
 
  The
 proof is the same as the proof for the first step in part \emph{i}) and will therefore be omitted. It follows that
 either $I=(p^m(b_0T_0 -2\overline{a}-\overline{c}T_1))$ or $I=0$. 
 
 \begin{bf} Claim:\end{bf} \emph{ $I \neq 0$.}
 
  Proceeding as in part \emph{i}) this will be done by constructing a lifting of the 
Hodge filtration of $x$ 
over $W[z]/(z^2-p,\ p^{m+1})$  which is stable under $ps_j$, but not 
stable under $j$.
Again let
$L=\mathbb{L}$ be the Dieudonn\a'e module of $x$. Then  a
basis of $L$ is given by $e_1, e_2, e_3= \Pi e_1,e_4= p^{-1}\Pi e_2$. Let 
\[ 
L_z=L\otimes_W W[z]/(z^2-p),
\] 
 and let 
 \[ 
P=L\otimes_W W[z]/(z^2-p,p^{m+1})=L_z/p^{m+1}L_z.
\] 
Since $\overline{c}$ and $b_0$ are units (see \cite{KR1}), we can choose $t \in W^{\times}$ such that 
$t^2 $ is not congruent to $\frac{\overline{c}}{b_0}$ modulo $p$.
Define 
\[
\begin{split} 
&f_0 \ = \ te_2  + z(e_1+e_2)\ \in L_z, \\
&f_1 \ = \ e_3  + pe_4 +zte_4 \ \in L_z.
\end{split}
\]
Denote again by $\overline{f_0}$ resp. $\overline{f_1}$ the image of $f_0$ resp. $f_1$ in $P$ and define
a Hodge filtration $\mathcal{F}\hookrightarrow P$ by $\mathcal{F} \ = \ < \overline{f_0}, \overline{f_1}>$. 
As above we see that it lifts the Hodge filtration of $x$, in particular that
 $\Pi f_0=zf_1$ and $\Pi f_1=zf_0$ holds. Using $p \mid \overline{a}$  and $ p\mid\overline{b}=pb_0$ 
 we
calculate in $L_z$:
\[
 ps_jf_0 \equiv p^mz\overline{c}t^{-1}f_0  \ \begin{rm}  mod \end{rm} \ p^{m+1} \ \text{ and } \ 
ps_jf_1 \equiv p^mz b_0 t f_1   \ \begin{rm}  mod \end{rm} \ p^{m+1},  
\]
showing that $ps_j\mathcal{F}\subset \mathcal{F}$. Writing $\overline{a}=pa_0$ we also calculate
\[
\Pi s_j f_1 \equiv p^mztb_0\cdot e_1 + (-p^mzta_0+p^m\overline{c})\cdot e_2  \ \begin{rm}  mod \end{rm} \ p^{m+1}.
\]
We claim now that $\Pi s_j \mathcal{F}\nsubseteq \mathcal{F}$. Otherwise we would have in particular 
$\Pi s_j\overline{f_1}\in
\mathcal{F}$. We then find $r \in W[z]/(z^2-p)$ with $ \Pi s_j f_1 \equiv rf_0 \begin{rm}  mod \end{rm} \
p^{m+1}$,
hence,
\[
\begin{split} 
r(t+z) & \equiv  -p^mzta_0+p^m\overline{c} \  \begin{rm}  mod \end{rm} \ p^{m+1}, \\
rz & \equiv p^mztb_0  \ \begin{rm}  mod \end{rm} \ p^{m+1}. 
\end{split}
\]
Subtracting the second congruence from the first and multiplying the result by $z$ we get
\[
rtz\equiv p^mz\overline{c} \ \begin{rm}  mod \end{rm} \ p^{m+1},
\]
 and on the other hand multiplying the second
congruence by $t$ we get 
\[
rtz\equiv p^mzt^2b_0  \ \begin{rm}  mod \end{rm} \ p^{m+1}.
\]
We conclude that $\overline{c}-t^2b_0\equiv 0\ \begin{rm}  mod \end{rm} \ z$ and hence that
$\overline{c}-t^2b_0\equiv 0\ \begin{rm}  mod \end{rm} \ p$, which contradicts the assumption we made on $t$.
Hence 
$\Pi s_j \mathcal{F} \nsubseteq \mathcal{F}$. This confirms the claim and ends the proof of \emph{ii}).
\newline

\emph{Proof} of \emph{iii}). In this case $\overline{a}$ is a unit and $\overline{b}$ is divisible by $p$, see \cite{KR1}. Since $\nu_p(\det(ps_j))\geq 1$, it follows that $m\geq 1$. We start as in part \emph{i}) by showing first that in
 some affine neighborhood the 
ideal $I$ of 
$Z(j)$ in $Z(ps_j)$ contains the ideal $p^m$. The proof is the same, so it will be omitted. 
Using the equations for $Z(ps_j)$, which are now given in a neighborhood of $x$ by $T_0p^m=T_1p^m=0$ 
(see \cite{KR1}),
 it follows that $I$
 is either the zero ideal or equals $(p^m)$. Thus we must show again that $I\neq 0$, which will  be done
 by constructing a lifting of the Hodge filtration of $x$ 
over $W[z]/(z^2-p,\ zp^{m+1})$  which is stable under $ps_j$, but not 
stable under $j$.
Again let
$L=\mathbb{L}$ be the Dieudonn\a'e module of $x$, and consider the
basis  $e_1, e_2, e_3= \Pi e_1,e_4= p^{-1}\Pi e_2$ of $L$. Let 
\[ 
L_z=L\otimes_W W[z]/(z^2-p),
\] 
 and let 
 \[ 
P=L\otimes_W W[z]/(z^2-p,zp^m)=L_z/zp^mL_z.
\] 
Define 
\[
\begin{split} 
&f_0 \ = \ e_2  + z(e_1+e_2)\ \in L_z, \\
&f_1 \ = \ e_3  + pe_4 +ze_4 \ \in L_z.
\end{split}
\]
Denote by $\overline{f_0}$ resp. $\overline{f_1}$ the image of $f_0$ resp. $f_1$ in $P$, and define
a  filtration $\mathcal{F}\hookrightarrow P$ by $\mathcal{F} \ = \ < \overline{f_0}, \overline{f_1}>$. 
By the same reasons as in part \emph{i})
it lifts the  Hodge filtration of $x$, and we have the relations $\Pi f_0=zf_1$ and $\Pi f_1=zf_0$.
By \cite{KR1}, in this case $\overline{b}$ is divisible by $p$ and $\overline{a}$ is a unit. We calculate
\[
 ps_jf_0 \equiv -p^m\overline{a}f_0  \ \begin{rm}  mod \end{rm} \ zp^m \ \text{ and } \ 
ps_jf_1 \equiv p^m\overline{a}f_1   \ \begin{rm}  mod \end{rm} \ zp^m,  
\]
showing that $ps_j\mathcal{F}\subset \mathcal{F}$. We also calculate,
\[
\Pi s_j f_1 \equiv p^m \overline{a}\cdot e_1 + (-p^{m-1}z\overline{a}+p^m\overline{c}-p^m\overline{a})\cdot e_2  \ \begin{rm}  mod
\end{rm} \ p^m.
\]
Proceeding as above we assume  $\Pi s_j \mathcal{F}\subset \mathcal{F}$ and have in particular 
$\Pi s_j\overline{f_1}\in
\mathcal{F}$. Hence we find $r \in W[z]/(z^2-p)$ with $ \Pi s_j f_1 \equiv rf_0 \  \begin{rm}  mod \end{rm} \
zp^m$ and hence
\[
\begin{split} 
(1+z)r & \equiv  -p^{m-1}z\overline{a}+p^m\overline{c}-p^m\overline{a} \  \begin{rm}  mod \end{rm} \ zp^{m}, \\
zr & \equiv p^m \overline{a} \ \begin{rm}  mod \end{rm} \ zp^{m}. 
\end{split}
\]
Subtracting the second congruence from the first
and multiplying the result by $z$, we get
\[
zr\equiv -p^m \overline{a} \ \begin{rm}  mod \end{rm} \ zp^m,
\]
which is a contradiction to the second congruence since $\overline{a}$ is a unit and $p \neq 2$. Hence 
$\Pi s_j \mathcal{F} \nsubseteq \mathcal{F}$.
This completes
the proof of the proposition.
\qed
\newline

Let $g\in \End(N,V)$. Following \cite{KR1} §4, we denote by
$Z(g)^{\begin{rm} pure \end{rm}}$  the closed subscheme of $Z(g)$ defined by the ideal sheaf 
of local
sections with finite support.

\begin{Cor}\label{div}
If $p > 3$, and $j$ is $*$-special with $j^2 \neq 0$ and $\nu_p(\begin{rm}  det \end{rm}(s_j)) \geq -1$, the antispecial  cycle $Z(j)$ is a divisor and equals $Z(ps_j)^{\begin{rm} pure  \end{rm}}$.
\end{Cor}
{\em Proof.} 
It follows from  the proof of Proposition \ref{ordloc}
 and from Proposition \ref{acycglss} 
that
\begin{equation}\label{puregl}
Z(j)^{\begin{rm} pure  \end{rm}}=
Z(j)=Z(ps_j)^{\begin{rm} pure  \end{rm}}.
\end{equation} \qed

%-----------------------------------------------------------------------------------------------------------------------
%-----------------------------------------------------------------------------------------------------------------------

\section{Intersection calculus of antispecial cycles}
We keep our fixed ramified automorphism $*$ of order $2$ of $B$.
In this section we calculate the intersection number of two antispecial cycles. 

On the space $\mathbb{Q}_p$-vector space $V[*]$ we have a quadratic form
 \[
 q(j)=(ps_j)^2= -p^2 \begin{rm}det\end{rm}(s_j).
 \]
($s_j$ is special and hence $s_j^2 \in \mathbb{Q}_p$, see \cite{KR1}, p. 167.)
Recall that $s_j=\iota(b_*^{-1})j$ and $b_*^2=\eta_*p$. 
We  also consider the quadratic form
 \[
 Q(j)=j^2=p^{-1}\eta_*q(j),
 \]

 Suppose we are given two  $*$-special endomorphisms $j_1$ and 
$j_2$. We assume that $j_i^2 \neq 0$ and $\nu_p(\begin{rm}  det \end{rm}((s_j)_i)) \geq -1$ for $i=1,2$ .
Following the general
definitions in \cite{KR1}, §4,
 we define the intersection number of their associated cycles by
\[
(Z(j_1),Z(j_2))=\chi(\mathcal{O}_{Z(j_1)}\otimes^{\mathbb{L}}\mathcal{O}_{Z(j_2)}),
\]
where $\chi$ denotes the \emph{Euler-Poincar\a'e} characteristic and $\otimes^{\mathbb{L}}$ denotes the derived tensor product.

By \cite{KR1}, Lemma 4.3, 
\begin{equation}\label{pure}
(Z(j_1),Z(j_2)) = (Z(j_1)^{\begin{rm} pure \end{rm}},Z(j_2)^{\begin{rm} pure \end{rm}}).
\end{equation}

We define $\bf{j}$ to be the $\mathbb{Z}_p$-span of $j_1$ and $j_2$ in $V[*]$ and we assume that
$\bf{j}$ has rank $2$.
Let $\beta$ be the bilinear form on $\bf{j}$ associated to the quadratic form $q$,
\[
\beta(x,y)=q(x+y)-q(x)-q(y).
\]
If $\bf{j}$ is nondegenerate with respect to this bilinear form, it follows from \cite{KR1}, Theorem 5.1 that
$(Z(j_1),Z(j_2))$ depends only on $\bf{j}$ (write $j_i=\iota(b_*)(s_j)_i$ and use that $(s_j)_i$ is special). 

On the other hand, since $p\neq 2$, we can choose a $\mathbb{Z}_p$-basis $j,j^{'}$ of $\bf{j}$, 
for which the matrix of $\beta$ 
is diagonal, 
\begin{align}
T:=
\begin{pmatrix}
 q(j)& \frac{1}{2}\beta(j,j^{'}) \\ 
 \frac{1}{2}\beta(j,j^{'}) &q(j^{'})
\end{pmatrix}
=\text{diag}(\eta_1p^{\beta_1},\eta_2p^{\beta_2})
\end{align}
where $\eta_1, \eta_2 \in 
\mathbb{Z}_p^\times$ and $\beta_1 \leq \beta_2$.

Define $\varepsilon_i\in \mathbb{Z}_p^\times $ and $\alpha_i \in \mathbb{N}$ by 
 $\eta_*\eta_i p^{\beta_i -1}=
\varepsilon_i p^{\alpha_i}$. Then $\alpha_1 \leq \alpha_2$. 
If $j_1,j_2$ is a basis of ${\bf j}$ for which $T$ has the form  $\text{diag}(\eta_1p^{\beta_1},\eta_2p^{\beta_2})$ then $Q(j_i)=j_i^2=\varepsilon_ip^{\alpha_i}.$ 
The numbers $\beta_i$ and hence the $\alpha_i$ are
invariants of $\bf{j}$. In case $\beta_1 < \beta_2$, the units $\eta_i$ and $\varepsilon_i$ are up to a square 
uniquely determined by
$\bf{j}$ as well.
Using (\ref{pure}) and (\ref{puregl}) we get from Theorem 6.1 of \cite{KR1} directly an expression for $(Z(j_1),Z(j_2))$ using the invariants 
$\beta_i$ and $\eta_i$ of $\bf{j}$ defined by the quadratic form $q$. For later use we translate this formula into an expression depending on the invariants $\alpha_i$ and $\varepsilon_i$ defined by the quadratic form $Q$.

\begin{The}\label{intthe}
Let $p>3$, let $j_1,j_2 \in V[*]$, and let $\bf{j}$ be their $\mathbb{Z}_p$-span in $\End(N,V)$. Assume that 
$\bf{j}$ 
is of dimension $2$ and
 nondegenerate. Using the same notations as above the following formula holds
\[
\begin{split}
&(Z(j_1),Z(j_2))= \\& \alpha_1+\alpha_2+3 
 -
\begin{cases}
p^{(\alpha_1 +1)/2}+2\frac{p^{(\alpha_1 +1)/2}-1}{p-1} \ & \text{if $\alpha_1$ is odd and 
$\chi(\eta_*\varepsilon_1)=-1$}\\
(\alpha_2-\alpha_1+1)p^{(\alpha_1 +1)/2}+ 2\frac{p^{(\alpha_1 +1)/2}-1}{p-1}& \text{if $\alpha_1$ is odd and 
$\chi(\eta_* \varepsilon_1)=1$}\\
2\frac{p^{\alpha_1/2+1}-1}{p-1} \ & \text{if $\alpha_1$ is even,}
\end{cases}
\end{split}
\]
where $\chi$ denotes the quadratic residue character on $\mathbb{Z}_p^\times$ resp.
$\mathbb{F}_p^\times$.
\end{The}

 (Note that in case $\alpha_1=\alpha_2$  this expression does not depend on $\varepsilon_1$. Note also that $\varepsilon_*$ is well defined up to multiplication by a unit of $\Zp$ and hence $\chi(\eta_*)$ is well defined.) 
 \qed
\begin{Rem}\label{pgd}{\em  Theorem \ref{intthe} is still valid for  $p=3$ provided that $\alpha_1 \geq 1$.
Indeed,
the only point where we used $p>3$ is in Lemma \ref{pdringe} for the
nilpotence of the $pd$-structure of the maximal ideals of 
rings of the form $A=W[x]/(x^2-p \varepsilon, x^r)$. These rings are only needed in the proof of Proposition
 \ref{acycglss}. Looking back at this proof  and using the same notation, 
 we see that in case (\emph{i}) and in case (\emph{iii}) we do not need these rings if we only want to
  determine 
 $Z(j)^{\begin{rm} pure  \end{rm}}$ in a neighborhood of the superspecial point $x$ we are considering.
  In case (\emph{ii}) rings of the above form are also needed to determine 
 $Z(j)^{\begin{rm} pure  \end{rm}}$ in a neighborhood of $x$ 
 if and only if $m=0$. But in this case we have $ \nu_p ((ps_j)^2)=1$, since 
 $\overline{a}$ is divisible by $p$ and $b_0$ and $c$ are units. Hence, if $p=3$ and $\alpha_1 \geq 1$ (and hence 
  $\beta_1 > 1$), we see that $Z(j_1)^{\begin{rm} pure  \end{rm}}$ and 
 $Z(j_2)^{\begin{rm} pure  \end{rm}}$ are given by the same equations as in case $p>3$.}
\end{Rem}

%-----------------------------------------------------------------------------------------------------------------------
%-----------------------------------------------------------------------------------------------------------------------

\section{An application to Arithmetic Hirzebruch-Zagier cycles}

We consider a supersingular  formal 
$p$-divisible group $\mathcal{A}$ over $k$ of height 4 and dimension 2 which is equipped with an action 
\[
\iota_0: \mathbb{Z}_{p^2}\rightarrow \begin{rm}End\end{rm}(\mathcal{A}),
\]
 such that $\mathcal{A}$ is special with respect to $\iota_0$. We further assume that $\mathcal{A}$ is equipped 
  with a  polarization

\begin{equation*}
\begin{CD}
\lambda: \mathcal{A} @ >\sim>> \hat{\mathcal{A}},
\end{CD}
\end{equation*}
such that for the Rosati involution $\iota_0(a)^{*}=\iota_0(a)$.

  We consider the following functor ${\mathcal{M}}^{HB}$ on the category 
$\begin{rm}Nilp \end{rm}$ of $W$-schemes $S$ such that $p$ is locally nilpotent in
$\mathcal{O}_S$. It associates to a scheme $S \in \Nilp$ 
 the set of isomorphism classes of the following data.

\begin{enumerate}[(1)]
\item A $p$-divisible group $X$ over $S$, with an action
\[
\iota_0: \mathbb{Z}_{p^2}\rightarrow \begin{rm}End\end{rm}(X),
\]
such that $X$ is special with respect to this $\mathbb{Z}_{p^2}$-action.
\item 
A quasi-isogeny of height zero
\[
\varrho: \mathcal{A}\times_{\Spec k}\overline{S} \rightarrow X\times_S \overline{S},
\]
which commutes with the action of $\mathbb{Z}_{p^2}$   such that the following condition holds.
Let   ${\lambda}_{\overline{S}}: \mathcal{A}_{\overline{S}} \rightarrow
 \hat{\mathcal{A}}_{\overline{S}}$ be the map induced by $\lambda$.
Then we require the existence of an isomorphism $\tilde{\lambda}
:X \rightarrow \hat{X}$
such that for the induced map 
  $\tilde{\lambda}_{\overline{S}}:X_{\overline{S}} \rightarrow \hat{X}_{\overline{S}}$ we have 
${\lambda}_{\overline{S}}=\hat{\varrho}\circ \tilde{\lambda}_{\overline{S}}\circ \varrho$.

\end{enumerate}

(Here, as in the preceding sections, a $p$-divisible group $X$ over $S$ with $\mathbb{Z}_{p^2}$-action is said to be special if the induced $\mathbb{Z}_{p^2}\otimes \mathcal{O}_S$-module $\Lie X$ is, locally on $S$, free of rank $1$.)

Denote
the isocrystal of $\mathcal{A}$ by $N$. From the      
polarization $\lambda$ we get we a perfect symplectic form $\langle  ,\rangle  $ on the Dieudonné module of $\mathcal{A}$ and hence also on $N$. 
Let $\Lambda_0$ be a $W$-lattice in $N$
 which is stable under $F$ and $V$ and under the action of $\mathbb{Z}_{p^2}$ and for which
the dual lattice equals $\Lambda_0$ (via the identification induced by $\langle  ,\rangle  $.)
Then setting $\mathcal{L}=\{\Lambda_0$\}  and $O_B=\mathbb{Z}_{p^2}$, the functor $\mathcal{M}^{HB}$ is a special case of \cite{RZ}, Definition 3.21. (In this definition there are imposed some additional conditions which are automatic here.)
By Theorem 3.25 of loc. cit. the functor ${\mathcal{M}}^{HB}$ is representable by a formal scheme which we
also call ${\mathcal{M}}^{HB}$. This formal scheme is formally locally of finite type over $W$ and is formally smooth  over $\mathbb{Z}_p$. 
%It has dimension $3$.

Following \cite{KR2}, §5, we define in this context the space of special endomorphisms
\[
V^{'}=\{j \in \End(N); j\iota_0(a)=\iota_0(a^{\sigma})j \ \text{ and } \ j^*=j\},
\]

 where $*$ denotes 
the adjoint with
respect to the alternating form $\langle  ,\rangle  $. As shown in loc. cit., $V^{'}$ is a $4$-dimensional vector space over
$\mathbb{Q}_p$ with  quadratic form 
\[
Q(j)=j^2.
\]

For $j \in V^{'}$  we define 
the special cycle $Z(j)$ associated to $j$ to be the closed formal subscheme of ${\mathcal{M}}^{HB}$ consisting of
all points $(X, \varrho)$ such that $\varrho \circ j \circ \varrho^{-1}$ lifts to an endomorphism of $X$.
Again, the fact that $Z(j)$ is a closed formal subscheme of ${\mathcal{M}}^{HB}$ follows from \cite{RZ}, Proposition 2.9.

We fix $j_1 \in V^{'}$  with $j_1^2=\varepsilon_1p$ for some unit
$\varepsilon_1\in \mathbb{Z}_p^{\times}$. (In \cite{KR2}, p. 188 the space $V^{'}$ is described more precisely
and from this description one easily sees that such $j_1$ exist.) Our next aim is to identify $Z(j_1)$ with the
Drinfeld moduli scheme $\mathcal{M}$ introduced in §2. Let $\varepsilon \in \mathbb{Z}_{p^2}^{\times}$ be such
that $\varepsilon \cdot \varepsilon^{\sigma}=\varepsilon_1$.
Let $(X,\varrho)\in  Z(j_1)$. We define an  $O_B$-operation $\iota$ on the points  $X$ by
keeping the $\mathbb{Z}_{p^2}$-action $\iota_0$ and by setting $\iota(\Pi)=\iota(\varepsilon^{-1})j_1$.  Since  $\mathcal{A}$ has height $4$, for  any point $(X,\varrho)$ the $p$-divisible group $X$ also has height $4$. Since $X$ is special, it has dimension $2$. 
 We must check that the condition given in (2) in the definition
of ${\mathcal{M}}^{HB}$ is automatic for $\mathcal{M}$. But this is done in the proof of Proposition 3.3,
Chapitre III of \cite{BC}. (We may suppose that $\mathbb{X}$ is superspecial. Then the diagram on p. 138 of loc. cit. in case $S=\Spec (B)$, where $B$ is a  $\Zp^{\rm nr}$-algebra in which  $p$ is nilpotent, is the  diagram that we need in this case (in loc. cit. the isomorphism $\tilde{\lambda}$ is called $p$). Since the solution for   $\tilde{\lambda}$ is unique in this case, we get the general case by glueing the solutions in formal neighborhoods of the geometric points.)

We define
\[
V^{'}[j_1]=\{j\in V^{'}\  | \ j \perp j_1 \text{ with respect to the bilinear form associated to }  Q \}.
\]

\begin{Lem}\label{spur0}
Let $j,j^{'} \in V^{'}$. Then $\iota_0(\delta)j^{'}j$ has trace $0$.
\end{Lem}
\emph{Proof.}
We can choose a basis $e_1,e_2,e_3,e_4$ of the 
Dieudonn\a'e module $M$ of $\mathcal{A}$ such that $e_1,e_2\in M_0$ and $e_3,e_4\in M_1$ and such that the matrix of $j$ resp. the matrix of $j^{'}$ has the form
\[
j=
\begin{pmatrix}
 & & d&-b \\
 & & -c&a \\
a&b & & \\
c&d & & 
\end{pmatrix},
 \text{ resp.  }
j^{'}=
\begin{pmatrix}
& & d^{'}&-b^{'} \\
& & -c^{'}&a^{'} \\
a^{'}&b^{'}& & \\
c^{'}&d^{'}& &
\end{pmatrix}.
\]
(Compare the proof of Proposition \ref{divisor} below or \cite{KR2}, §5.) The diagonal entries of $jj^{'}$ are $(da^{'}-bc^{'}),(-cb^{'}+ ad^{'}),( ad^{'}-bc^{'}),(-cb^{'}+da^{'}).$
Since $\iota_0(\delta)$ acts on $M_0$ by    
multiplication by $\delta$ and on $M_1$ by multiplication by $-\delta$, the claim follows.
\qed
\newline

\begin{Pro}\label{idefix}
  Let $*=\Int(\delta j_1)$ where we identify $B$ with  $ \mathbb{Q}_{p^2}[j_1]$. If  $j \in V^{'}[j_1]$ then $j$ is $*$-special. 
\end{Pro}
\emph{Proof.} Since for $j\in V^{'}[j_1]$ we have $j\iota(\delta)j^{-1}=\iota(-\delta)$ and $jj_1j^{-1}=-j_1$, it follows that $j\iota(a)j^{-1}=\iota(a^*)$ for all $a\in B$. Using lemma \ref{spur0} the claim follows. \qed
\newline

From this proposition we conclude

\begin{Cor}\label{idefix1}
Let  $j \in V^{'}[j_1]$. Identifying $Z(j_1)$ with $\mathcal{M}$ as above, the intersection $Z(j_1)\cap Z(j)$  equals the antispecial cycle $Z(j)$ in $\mathcal{M}$. \qed
\end{Cor}

\begin{Rem}
\emph{In \cite{KR2}, p. 243  it is erroneously asserted that conjugation by $j \in V^{'}[j_1]$ induces the main involution on $B$.}
\end{Rem}

\begin{Pro}\label{divisor}
 Let $j \in V^{'}$
 be such that $Q(j)\neq0$ and $Z(j)\neq \emptyset$. 
 Then
   $Z(j)$ is a divisor in ${\mathcal{M}}^{HB}$.
\end{Pro}
\emph{Proof}. 
Let $x\in {\mathcal{M}}^{HB}$ be a closed point which belongs to $Z(j)$, and
let $R$ be the local ring $\mathcal{O}_{{\mathcal{M}}^{HB},x}$. Let  $J$ be the ideal of $R$ coming from the
ideal sheaf of $Z(j)$, and let $\mathfrak{m}$ be the maximal ideal of $R$. We must show that $J$ is a principal ideal. 
Let $A:=R/(\mathfrak{m}J)$, and let $\overline{A}:=R/J$. We have $\overline{A}=A/I$, where $I={J}/({\mathfrak{m}}{J})$.
By Nakayama's lemma it is enough to show that $I$ is a principal ideal.
 Since $I^2=0$, the ideal $I$ carries a nilpotent $pd$-structure. 
Now  $A$ is separable and complete for the topology defined by the ideal $(p)$. Therefore,  by considering projective limits, it follows that we can apply Grothendieck-Messing theory for the pair $A$, $\overline{A}$. 
There is an $A$-valued and  an $\overline{A}$-valued point of ${\mathcal{M}}^{HB}$ in the natural way.   The latter also gives an $\overline{A}$-valued point of $Z(j)$. 
Let $M$ be the value in $A$ of the crystal of the $p$-divisible group belonging to the  $A$-valued point of ${\mathcal{M}}^{HB}$, and let $\overline{M}$ be the value in $\overline{A}$ of the  crystal of the $p$-divisible group belonging to the  $\overline{A}$-valued point of ${\mathcal{M}}^{HB}$.
These are free modules of rank $4$ over $A$ resp. $\overline{A}$, which are equipped with a perfect alternating form $\langle  ,\rangle  $ and a $\mathbb{Z}_{p^2}$-action $\iota_0$ such that $\iota_0(a)$ is selfadjoint. We have $\overline{M}=M \otimes_{A}\overline{A}$. From the $\mathbb{Z}_{p^2}$-actions we get $\mathbb{Z}/2$-gradings $M=M_0\oplus M_1$ resp. $\overline{M}=\overline{M}_0\oplus \overline{M}_1$.
 
  Let $\overline{\mathcal{F}}\hookrightarrow \overline{M}$ 
be the Hodge filtration over $\overline{A}$ corresponding to the $\overline{A}$-valued point of $Z(j)$. 
(The submodule $\overline{\mathcal{F}}$ is free of rank $2$ and stable under $\mathbb{Z}_{p^2}$ and under $j$.)   By Grothendieck-Messing Theory, to a lifting of the corresponding $p$-divisible group over $A$
corresponds a lifting of the Hodge filtration  over $A$.
Let ${\mathcal{F}}\hookrightarrow {M}$ be the Hodge filtration corresponding to the $A$-valued point of ${\mathcal{M}}^{HB}$. This lifts the Hodge filtration over $\overline{A}$. 

Let $\{e_1,e_2,e_3,e_4\}$ 
be a basis of $M$ such that $e_1,e_2 \in M_{0}$ and $e_3,e_4\in M_{1}$.
Then the images $\overline{e_i}$ of the $e_i$ in $\overline{M}$  form a basis of
$\overline{M}$. We may suppose  that $\{\overline{e_2},\overline{e_3}\}$ is a basis of $\overline{\mathcal{F}}$.
Since $M_{0} \perp M_{1}$ (which follows from the fact that $\iota_0(\delta)$ is selfadjoint), and since the determinant of the matrix
 of the bilinear form $\langle  ,\rangle  $ is a unit, we may
suppose that $\langle  e_1,e_2\rangle  =\langle  e_3,e_4\rangle  =1$. As in \cite{KR2}, §5 it follows that with respect to the basis $e_1,...,e_4$ the matrix of $j$ has the form

\[
j= \begin{pmatrix}
 &j_1 \\
  j_0&
\end{pmatrix},
\]
where $j_i \in \begin{rm} Hom\end{rm}(M_{i},M_{i+1})$ and $j_1=j_0^*$. Further for
\[
j_0=
\begin{pmatrix}
 a&b \\
  c&d
\end{pmatrix},
\]
we have 
\[
j_1=j_0^*=
\begin{pmatrix}
 d&-b \\
  -c&a
\end{pmatrix}.
\]

Since $\mathcal{F}$ lifts $\overline{\mathcal{F}}$  we find a basis 
$f_0,f_1$ of ${\mathcal{F}}$ such that $f_0=e_2+m_1e_1$ and $f_1=e_3+m_4e_4$ for some $m_1,m_4 \in I$.

Now, let $\mathfrak{b}\subset I$ be an ideal in $A$ and let $B=A/ \mathfrak{b}$.

\begin{bf} Claim: \end{bf} \emph{The map $\Spec B \rightarrow \mathcal{M}^{HB}$ factors through $Z(j)$ if and only if $bm_4-cm_1-d=0$ in $B$.}

By the same reasons as above we may apply Grothendieck-Messing theory to the pairs $B, \overline{A}$ and $ A,B$.
Let $M_B$ be the value in $B$ of the crystal of the $p$-divisible group belonging to the  $B$-valued point of ${\mathcal{M}}^{HB}$. 
 It equals $M\otimes_A B$. To the lifting $\Spec B \rightarrow \mathcal{M}^{HB}$ 
of $\Spec \overline{A} \rightarrow \mathcal{M}^{HB}$ corresponds the Hodge    
 filtration $\mathcal{F}_B \hookrightarrow M_B$,
 where $\mathcal{F}_B=\mathcal{F} \otimes_A B$. (Note that  the map
 $\Spec A \rightarrow \mathcal{M}^{HB}$ lifts $\Spec B \rightarrow \mathcal{M}^{HB}$.)

In $M_B$ we have $j(f_0)=be_3+de_4+m_1(ae_3+ce_4)$. 
The image  $j(f_0)$ 
 lies in $\mathcal{F}_B$ if and only if $j(f_0)=xf_1$ in $\mathcal{M}_B$ for some
$x\in B$. That means $be_3+de_4+m_1(ae_3+ce_4)=x(e_3+m_4e_4)$ in $\mathcal{M}_B$, hence $x=b+am_1$ and $xm_4=d+cm_1$ in $B$, and hence 
\[
j(f_0)\in \mathcal{F}\Leftrightarrow bm_4-cm_1-d=0 \text{ in } B.
\] 
 Evaluating the corresponding condition for
$j(f_1)$ one gets the same equation.
This confirms the claim.

In case $\mathfrak{b}=I$ it follows in particular that $ bm_4-cm_1-d=0$ in $\overline{A}$.  Hence  in $A$ we have an inclusion $(bm_4-cm_1-d) \subset I $.
Therefore we can apply the  claim in case $\mathfrak{b}=(bm_4-cm_1-d)$ and conclude that the map  $\Spec A/(bm_4-cm_1-d)\rightarrow  \mathcal{M}^{HB}$ factors through $Z(j)$. 
On the other hand, by definition, $I$  is the minimal ideal of $A$ with the property that $\Spec A/I \rightarrow \mathcal{M}^{HB}$ factors through $Z(j)$. It follows that $I=(bm_4-cm_1-d)$.

It remains to show that the equation for $Z(j)$ is nowhere trivial. If the equation is trivial in some local ring of $\mathcal{M}^{HB}$ then it follows that there is some (open) formal affine neighborhood in which the equation for $Z(j)$ is trivial. If the equation for $Z(j)$ is not trivial in some local ring of $\mathcal{M}^{HB}$ then it follows that there is some (open) formal affine neighborhood in which  $Z(j)$ is given by one non trivial equation. It follows that the set of points of $\mathcal{M}^{HB}$ in whose local ring the equation for $Z(j)$ is trivial and the set of points of $\mathcal{M}^{HB}$ in whose local ring the equation for $Z(j)$ is not trivial are both open. By \cite{KR2}, Lemma 8.2 the latter set is not empty. Since $\mathcal{M}^{HB}$ is connected, it follows that the equation for $Z(j)$ is nowhere trivial.
\qed

\begin{Lem}\label{Lweg}
Suppose that $j \in V^{'}[j_1]$. Then 
\[
\mathcal{O}_{Z(j_1)} \otimes_{\mathcal{O}_{\mathcal{M}^{HB}}}^{\mathbb{L}} \mathcal{O}_{Z(j)}= \mathcal{O}_{Z(j_1)} \otimes_{\mathcal{O}_{\mathcal{M}^{HB}}}
\mathcal{O}_{Z(j)}.
\]
 More precisely, the object of the derived category on the l.h.s. is represented by the r.h.s..
\end{Lem}

\emph{Proof.} By Proposition \ref{divisor} we know that $Z(j_1)=\mathcal{M}$ is a divisor in $\mathcal{M}^{HB}$. Hence, locally on  $\mathcal{M}^{HB}$, 
there is an exact sequence
\begin{equation*}
\begin{CD}
0 @ >>> \mathcal{O}_{{\mathcal{M}}^{HB}} @>f>>\mathcal{O}_{{\mathcal{M}}^{HB}} @>>>
\mathcal{O}_{{\mathcal{M}}}@>>> 0,
\end{CD}
\end{equation*}
where $f=0$ is the equation of $Z(j_1)$. Using this resolution we see that 
$\mathcal{O}_{Z(j_1)} \otimes_{\mathcal{O}_{\mathcal{M}^{HB}}}^{\mathbb{L}} \mathcal{O}_{Z(j)}$ is represented by the complex 
$
\begin{CD}
\mathcal{O}_{Z(j)} @>f>>\mathcal{O}_{Z(j)}.
\end{CD}
$
Its cohomology sheaves are $\mathcal{O}_{Z(j)} \otimes \mathcal{O}_{Z(j_1)}$ and $\mathcal{T}or_1^{\mathcal{O}_{{\mathcal{M}}^{HB}}}(\mathcal{O}_{Z(j_1)}, \mathcal{O}_{Z(j)})$. 
 
 Hence we only have to show that
 $\mathcal{T}or_1^{\mathcal{O}_{{\mathcal{M}}^{HB}}}(\mathcal{O}_{Z(j_1)},\mathcal{O}_{Z(j)}) =0$. 
 So we have to
 show that in every local ring of $Z(j)$ the image of $f$ is not a zero divisor. For this we consider a local ring of ${\mathcal{M}}^{HB}$. Let $g=0$ be the equation for $Z(j)$ in this local ring. (The proof of Proposition \ref{divisor} shows that $Z(j)$ is in fact given by one equation in any local ring of ${\mathcal{M}}^{HB}$.) Then we have to show that $g$ and the image of $f$ in the local ring have no common prime divisor.
 (Since ${\mathcal{M}}^{HB}$ is regular, its local rings are unique factorization domains.) Assuming the
 contrary it would follow that $\mathcal{M}$ and $Z(j)$ have a common component. But this
 contradicts the fact that their intersection is pure one dimensional by the results of section $2$, which we can apply because of Corollary \ref{idefix1}.
 \qed
\newline

Let $j_2,j_3\in V^{'}[j_1]$ be such that the $\mathbb{Z}_p$-span $\begin{bf} j \end{bf}=
\mathbb{Z}_pj_2 + \mathbb{Z}_pj_3$ has rank $2$ as a submodule of $V^{'}$ and such that $Q$ induces a nondegenerate bilinear form on $\bf{j}$. 
We further suppose
that the matrix of the bilinear form on $\bf{j}$ associated
to $Q$ with respect to the basis $  j_2, \ j_3$ is equivalent to 
 $\begin{rm}  diag \end{rm}( \varepsilon_2 p^{\beta_2}, \varepsilon_3 p^{\beta_2})$ with $\varepsilon_i\in \Zp^{\times}$ and $1\leq \beta_2\leq \beta_3$.
  In this situation we define the intersection product of $Z(j_1)$, $Z(j_2)$ and $Z(j_3)$ by the \emph{Euler-Poincar\a'e} characteristic of the derived tensor product,
\[
(Z(j_1),Z(j_2),Z(j_3)):= \chi(\mathcal{M}^{HB},\mathcal{O}_{Z(j_1)}\otimes^{\mathbb{L}}\mathcal{O}_{Z(j_2)}
\otimes^{\mathbb{L}}\mathcal{O}_{Z(j_3)}).
\]\
This is well defined since $Z(j_1)\cap Z(j_2)\cap Z(j_3)$ is proper over $\Spec k$. Indeed, $Z(j_1)\cap Z(j_2)\cap Z(j_3)=(\mathcal{M}\cap Z(j_2))\cap (\mathcal{M}\cap Z(j_3))$).
This is included in  $ Z(\Pi j_2)\cap  Z(\Pi j_3)$ regarded as an intersection of special cycles inside $\mathcal{M}$. But this is proper over $\Spec k$ by \cite{KR1}, Proposition 3.6 and Corollary 2.14.

\begin{Pro}\label{vergl}
There is an equality of intersection multiplicities on $\mathcal{M}^{HB}$ resp. $\mathcal{M}$,
\[
(Z(j_1),Z(j_2),Z(j_3))=
((\mathcal{M}\cap Z(j_2)),(\mathcal{M}\cap Z(j_3))),
\]
where we regard the intersections $\mathcal{M}\cap Z(j_i)$ as antispecial cycles in $\mathcal{M}$, cf. Corollary \ref{idefix1}.

{\em The latter intersection multiplicity is given explicitly by Theorem \ref{intthe} for $*=\Int(\delta j_1)$ after replacing $j_1, j_2$ in
  this theorem by $j_2,j_3$. }
\end{Pro}
\emph{Proof.}
We have  
\[
\begin{split}
& \chi(\mathcal{M}^{HB},\mathcal{O}_{Z(j_1)}\otimes^{\mathbb{L}}\mathcal{O}_{Z(j_2)}
\otimes^{\mathbb{L}}\mathcal{O}_{Z(j_3)})  \\
= & \chi(\mathcal{M}^{HB},(\mathcal{O}_{Z(j_1)}\otimes^{\mathbb{L}}\mathcal{O}_{Z(j_2)})
\otimes^{\mathbb{L}}_{\mathcal{O}_{Z(j_1)}}({\mathcal{O}_{Z(j_1)}}\otimes^{\mathbb{L}}\mathcal{O}_{Z(j_3)})).
\end{split}
\]
By Lemma \ref{Lweg}, and since $Z(j_1)=\mathcal{M}$, the latter expression equals  
\[
\begin{split}
& \chi(\mathcal{M},(\mathcal{O}_{\mathcal{M}}\otimes \mathcal{O}_{Z(j_2)})\otimes_{\mathcal{O}_{\mathcal{M}}}
^{\mathbb{L}} (\mathcal{O}_{\mathcal{M}}\otimes \mathcal{O}_{Z(j_3)}))=
\chi(\mathcal{M}, \mathcal{O}_{\mathcal{M}\cap Z(j_2)}\otimes^{\mathbb{L}}\mathcal{O}_{\mathcal{M}\cap Z(j_3)}) \\
=&((\mathcal{M}\cap Z(j_2)),(\mathcal{M}\cap Z(j_3))).
\end{split}
\]
Using Proposition \ref{idefix} and Corollary \ref{idefix1} it follows that
 this
 intersection multiplicity can be calculated as in Theorem \ref{intthe} for $*=\Int(\delta j_1)$.
\qed

%-----------------------------------------------------------------------------------------------------------------------
%-----------------------------------------------------------------------------------------------------------------------

\section{Representation densities}
We recall that, for $S\in \begin{rm}  Sym \end{rm}_m(\mathbb{Z}_p)$ and 
$T\in \begin{rm}  Sym \end{rm}_n(\mathbb{Z}_p)$
with $\begin{rm}  det \end{rm}(S)\neq 0$ and $\begin{rm}  det \end{rm}(T)\neq 0$,
 the representation density is defined as 
\[
\alpha_p(S,T)= 
\operatorname*{lim}_{t\rightarrow\infty} p^{-tn(2m-n-1)/2} \mid \{x \in M_{m,n}(\mathbb{Z}/p^t\mathbb{Z});
\ S[x]-T \in p^t\begin{rm}  Sym \end{rm}_n(\mathbb{Z}_p)\}\mid.
\]
Given $S$ as above, let 
\[
S_r=
\begin{pmatrix}
 S \\ 
 & 1_r \\
 & & -1_r
\end{pmatrix}.
\]
Then there is a rational function $A_{S,T}(X)$ of $X$ such that 
\[
\alpha_p(S_r,T)= A_{S,T}(p^{-r}).
\]
Let 
\[
\alpha_p^{'}(S,T)=\frac{\partial}{\partial X}(A_{S,T}(X)) \arrowvert_{X=1}.
\]
(Comp. \cite{KR1}, §7.)

Let $j_1$ be as in section $4$, and let $j_2,j_3 \in V^{'}[j_1]$  also be 
as in section $4$, i.e. 
such that the $\mathbb{Z}_p$-span 
$\begin{bf} j \end{bf}=  \mathbb{Z}_p  j_2 + \mathbb{Z}_p  j_3$  is 
 of rank $2$  and
nondegenerate for the bilinear form associated to the quadratic form $Q$, and such that the matrix of this bilinear form on $\bf{j}$ with respect to the basis $  j_2, \ j_3$ is equivalent to 
 $\begin{rm}  diag \end{rm}( \varepsilon_2 p^{\beta_2}, \varepsilon_3 p^{\beta_2})$ with $\varepsilon_i\in \Zp^{\times}$ and $1\leq \beta_2\leq \beta_3$. 
 Let $S=\begin{rm}  diag \end{rm}(1,-1,1,-\Delta)$. Then $T=\begin{rm}  diag \end{rm}( \varepsilon_1 p, \varepsilon_2 p^{\beta_2}, \varepsilon_3 p^{\beta_2})$ is represented by the space $V^{'}$, hence
 $T$ is not represented  by  $S$, see \cite{Ku}, Proposition 1.3. 

\begin{The}\label{theo}
Using the notation just introduced we have
\[
(Z(j_1),Z(j_2),Z(j_3))= -\frac{p^4}{(p^2+1)(p^2-1)}\alpha_p^{'}(S,T).
\]
\end{The}
\emph{Proof.} Let $\eta\in \mathbb{Z}_p^{\times}$, and let $S(\eta)=\begin{rm}  diag \end{rm}(1,-1,1,-\eta)$ so that $S=S(\Delta)$. 
Then it follows from \cite{S}, Lemma
3.5, that there is a polynomial $g_{T}(X)\in \mathbb{Z}[X]$ such that
$
\alpha_p(S(\eta)_r,T)=g_T(\chi(\eta)\cdot p^{-r-2}).
$
On the other hand, by \cite{Ka} there exists a polynomial $f_T(X)\in \mathbb{Q}[X]$ with 
$\alpha_p(S(1)_r,T)=f_T(p^{-r})$. Hence for any $r\in \mathbb{N}$ we have $g_T(p^{-r-2})=f_T(p^{-r})$ and hence
$g_T(p^{-2}X)=f_T(X)$. Therefore $\alpha_p(S_r,T)=f_T(-p^{-r})$ and hence $A_{S,T}(X)=f_T(-X)$.
 Katsurada's polynomial $f_T(X)$ is given explicitly in
\cite{W}. Following this article, we can write 
\[
f_T(X)=\tilde{\gamma}_p(T;X)\tilde{F}_p(T;X),
\]
where 
\[
\tilde{\gamma}_p(T;X)=(1-p^{-2}X)(1-p^{-2}X^2)
\]
and where 
$\tilde{F}_p(T;X)$ is defined as follows. First we define some invariants of
 $T=\begin{rm}  diag \end{rm}(\varepsilon_1 p, \varepsilon_2 p^{\beta_2}, \varepsilon_3 p^{\beta_2})$. Let
 \[
 \tilde{\xi}=
 \begin{cases}
\chi(-\varepsilon_1\varepsilon_2)\ & \text{if $\beta_2$ is odd,} \\
0 \ & \text{if $\beta_2$ is even,} 
 \end{cases}
 \]
and let 
\[
\sigma=
 \begin{cases}
2\ & \text{if $\beta_2$ is odd,} \\
1 \ & \text{if $\beta_2$ is even.} 
 \end{cases}
 \]
Further, let 
\[
\eta=
 \begin{cases}
+1\ & \text{if $T$ is isotropic,} \\
-1 \ & \text{if $T$ is anisotropic.} 
 \end{cases}
 \]
By \cite{W}, 2.11, we then have
\[
\begin{split}
\tilde{F}_p(T;X) &= \sum_{i=0}^{1}\sum_{j=0}^{(1+\beta_2-\sigma)/2-i}p^{i+j}X^{i+2j} \\
& + \eta \sum_{i=0}^{1}\sum_{j=0}^{(1+\beta_2-\sigma)/2-i}p^{(1+\beta_2-\sigma)/2-j}X^{\beta_3+\sigma+i+2j} \\
& +  \tilde{\xi}^2p^{(1+\beta_2-\sigma+2)/2}\sum_{i=0}^{1}\sum_{j=0}^{\beta_3-\beta_2+2\sigma-4}\tilde{\xi}^j
X^{\beta_2-\sigma+2+i+j}.
\end{split}
\]
To distinguish whether $T$ is isotropic or anisotropic we recall the following fact (see \cite{W}, p. 189). Let $i,j \in \{1,2,3 \}$
with $i \neq j$ and $\beta_i \equiv \beta_j \mod 2$, and define $k\in\{1,2,3 \}$ by $\{i,j,k\}=\{1,2,3\}$. Then
$T$ is isotropic if and only if $\chi(-\varepsilon_i\varepsilon_j)=1$ or $\beta_k \equiv \beta_j \mod 2$. On
the other hand, since $T$ is represented by $V^{'}$ we have
\begin{align}\label{iden}
-1=(-1)^{1+\beta_2+\beta_3}\chi(-1)^{1+\beta_2+\beta_3+1\cdot\beta_2+1\cdot\beta_3+\beta_2\beta_3}
\chi(\varepsilon_1)^{\beta_2+\beta_3}\chi(\varepsilon_2)^{1+\beta_3}\chi(\varepsilon_3)^{1+\beta_2},
\end{align}
see \cite{Ku}, (1.16).

Now we evaluate 
$\alpha_p^{'}(S,T)=\frac{\partial}{\partial X}(\tilde{\gamma}_p(T;-X)\tilde{F}_p(T;-X))\arrowvert_{X=1}$. We show that
 $\alpha_p^{'}(S,T)$ equals $-(1+p^{-2})(1-p^{-2})$ times the expression given in   expression given in Theorem \ref{intthe} in case $*=\Int(\delta j_1)$ (comp. Proposition 
\ref{vergl})
after replacing $\varepsilon_1$ resp. $\varepsilon_2$ and $\alpha_1$ resp. $\alpha_2$ in
Theorem \ref{intthe} by $\varepsilon_2$ resp. $\varepsilon_3$ and $\beta_2$ resp. $\beta_3$. In the notation of Theorem \ref{intthe} we have $\eta_* = -\Delta \varepsilon_1$ . (Note that by Remark \ref{pgd} (using $\beta \geq 1$) we do not need to exclude the case $p=3$.)
Using these substitutions we distinguish the same cases as in Theorem \ref{intthe}.

\begin{bf} First case:\end{bf}\emph{ $\beta_2$ is odd and 
$\chi(-\Delta \varepsilon_1\varepsilon_2)=-1$, i.e. 
$\chi(-\varepsilon_1 \varepsilon_2)=1$.}

It follows immediately that $\tilde{\xi}=1$ and $\sigma=2$.

By (\ref{iden})  we get $-1=(-1)^{\beta_3}\chi(-\varepsilon_1\varepsilon_2)^{1+\beta_3}$,
 hence $\beta_3$ is odd.
>>From the criterion above we see that $T$ is isotropic and hence $\eta=1$. Now an easy calculation shows 
$\tilde{F}_p(T;-1)=0$ and 
\[
\frac{\partial}{\partial X}\tilde{F}_p(T;-X)\arrowvert_{X=1}=
-\beta_2-\beta_3-3 + p^{(\beta_2 +1)/2}+2\frac{p^{(\beta_2 +1)/2}-1}{p-1},
\]
which yields the claim in this case.

\begin{bf} Second case:\end{bf} \emph{ $\beta_2$ is odd and 
$\chi(-\Delta \varepsilon_1\varepsilon_2)=1$, i.e. 
$\chi(-\varepsilon_1 \varepsilon_2)=-1$.}

It follows immediately that $\tilde{\xi}=-1$ and $\sigma=2$.
We distinguish the subcases {\it $\beta_3$ even} and  {\it $\beta_3$ odd}.

\emph{If $\beta_3$ is even} we see 
 from the criterion 
above  that $T$ is anisotropic and hence $\eta=-1$.

\emph{If $\beta_3$ is odd} the criterion shows that $T$ is isotropic and hence $\eta=1$.

In both cases an easy calculation shows 
$\tilde{F}_p(T;-1)=0$ and 
\[
\frac{\partial}{\partial X}\tilde{F}_p(T;-X)\arrowvert_{X=1}=
-\beta_2-\beta_3-3 + (\beta_3-\beta_2+1)p^{(\beta_2 +1)/2}+ 2\frac{p^{(\beta +1)/2}-1}{p-1},
\]
which yields the claim in this case.

\begin{bf} Third case:\end{bf} \emph{$\beta_2$ is even.}
 
It follows immediately that $\tilde{\xi}=0$ and $\sigma=1$.
We distinguish the subcases {\it $\beta_3$ even} and { \it $\beta_3$ odd}.

\emph{If $\beta_3$ is even} we get 
>>from (\ref{iden}) that
$\chi(-\varepsilon_2\varepsilon_3)=1$, and by the criterion above $T$ is isotropic in this case
 and hence $\eta=1$. 

\emph{If $\beta_3$ is odd} we get 
>>from (\ref{iden}) that
$\chi(-\varepsilon_1\varepsilon_3)=-1$, and by the criterion above $T$ is anisotropic in this case and hence
$\eta=-1$.

In both cases an easy calculation shows 
$\tilde{F}_p(T;-1)=0$ and 
\[
\frac{\partial}{\partial X}\tilde{F}_p(T;-X)\arrowvert_{X=1}=
-\beta_2-\beta_3-3 + 2\frac{p^{\beta_2/2+1}-1}{p-1},
\]
which yields the claim in this case.
\qed

\end{document}